# Probabilistic and asymptotic methods with the Perron Frobenius's operator

Guy Cirier
guy.cirier@gmail.com

**Résumé**

*1- Écart résolvant. Soit $f$ une application polynomiale d'un compact $C$ de $\mathbb{R}^d$ dans lui-même. Une mesure $P$ est dite invariante sous $f$ ssi $P = P_f$ où $P_f = P \circ f^{-1}$. $P_f$ est appelé l'opérateur de Perron Frobenius (PF). Soient $e_f^n(y) = \partial^n \left[ e^{ya} - e^{yf(a)} \right] / \partial a^n \big|_{a=0}$ et $H_n(y) = \partial^n e^{yf(a)} / \partial a^n \big|_{a=0}$ où $0$ est point fixe de $f$. Soit $E_n$ l'idéal engendré par $d+1$ polynômes $e_f^n(y) = 0$ et $e_f^{n+1_1}(y) = 0$. Sous cette hypothèse et via la transformation de Laplace-Steiltjes, la solution est obtenue à partir de la distribution asymptotique des zéros réels de $E_n$ quand $n \to \infty$. En effet, l'idéal $E_n$ définit des variétés algébriques affines réelles. Si ces variétés sont $0$-dimensionnelles et ont une distribution asymptotique de densité dérivable $q$ $d$-dimensionnelle, la densité invariante sera de la forme:*

$$p(a) = \prod_{\square=1}^{\square=d} (-a_\square) \partial^d q(a) / \partial a_1 .. \partial a_\square .. \partial a_d$$

*2-Approximation. Sous les hypothèses de la méthode du col, la distribution asymptotique des zéros des $H_n(y)$ (variétés algébriques affines) est image réciproque de distributions uniformes $\kappa_\square$ sur $(0,1)$ :*

$$\pi \kappa_1 = s_1 \operatorname{Im} f_1(a) - \vartheta_1 \quad \square = 1, 2, ..., p$$

*pour toutes les $p \leq d$ coordonnées complexes du point critique $a$ de $\gamma(a) = s\zeta f(a) - \zeta \log a$. $\vartheta_\square$ est l'argument de $a_1$, $s_1 = y_1 / n$). Si la hessienne de $\gamma(a)$ est dégénérée, comme dans le cas de linéarités partielles, nous obtenons des familles de variétés aléatoires (exemple Henon).*

*3- Applications aux EDO ou EDP de la forme $\partial a / \partial t = F(a)$. On associe l'itération infinitésimale $f(a, \delta) = a + \delta F(a)$ pour lui appliquer l'opérateur de Perron Frobenius. Si $F(a)$ est partiellement linéaire, alors une distribution aléatoire peut être solution asymptotique de cette équation. Les applications les plus remarquables de l'étude sont les solutions asymptotiques des équations de Lorenz, de Navier Stokes ou d'Hamilton.*

*4- Résonance. Les valeurs propres $\lambda$ de la partie linéaire de $f$ vérifient $\lambda^k = 1$. Sous les hypothèses précédentes, le point critique $a$ de $\gamma(a)$ vérifie alors la condition of résonance*

$$1 = a^{k-1}$$

**Abstract**

*1- Solving deviation. Let $f$ be a polynomial application of a compact $C$ of $\mathbb{R}^d$ in $C$. A measure $P$ is invariant under $f$ iff $P = P_f$ where $P_f = P \circ f^{-1}$. $P_f$ is said Perron Frobenius operator (PF). Denote $e_f^n(y) = \partial^n \left[ e^{ya} - e^{yf(a)} \right] / \partial a^n \big|_{a=0}$ and $H_n(y) = \partial^n e^{yf(a)} / \partial a^n \big|_{a=0}$ where $0$ is a fixed point of $f$. Let $E_n$ be the ideal spanned by the $d+1$ polynomials $e_f^n(y) = 0$ and $e_f^{n+1_1}(y) = 0$. Under this hypothesis and via the Laplace-Steiltjes transform, we prove that the solution is given by the asymptotic distribution of the zeros of $E_n$ when $n \to \infty$. The ideal $E_n$ defines real algebraic affine manifolds. If these manifolds are 0-dimensional and have an asymptotic distribution with a derivable $d$-dimensional density $q$, the invariant density is:*

$$p(a) = \prod_{\square=1}^{\square=d} (-a_\square) \partial^d q(a) / \partial a_1 .. \partial a_\square .. \partial a_d$$

*2-Approximation. Under assumptions of the steepest descent, the asymptotic distribution of the zeros of the $H_n(y)$ (algebraic affine manifolds) getting the solution of (PF), is a function of iid $\kappa_l$ on (0,1)*

$$\pi \kappa_l = s_l \operatorname{Im} f_l(a) - \vartheta_l, \quad l = 1, 2, ..., p$$

*for all the $p$ complex coordinates of the critical point $a$ of $\gamma(a) = s\zeta f(a) - \zeta \log a$. $\vartheta_l$ is the argument of $a_l$, $s_l = \lim_{N \to \infty} y_l / n$. If the Hessian is degenerated, which is the case of partial linearity, we get a family of random manifolds as in the Henson's case.*

*3 – Applications to EDP, EDO. Asymptotic behaviours of ODE or PDE, as $\partial a / \partial t = F(a)$, via the Perron Frobenius's density, are most interesting. Here, the infinitesimal iteration is $f(a, \delta) = a + \delta F(a)$. But, if the Hessian is degenerated, the asymptotic design can be probabilistic. Among applications, are asymptotic solutions of Lorenz, Navier-Stokes or Hamilton's equations.*

*4- Resonance. If the eigenvalues $\lambda$ of the linear part of $f$ verify $\lambda^k = 1$. Under the previous hypothesis when $n = kq \to \infty$, at the critical point $a$ of $\gamma(a)$, we have the condition of resonance:*

$$1 = a^{k-1}$$

**Introduction**

The Perron-Frobenius's equation (PF) is a very difficult functional equation. Ulam and Von Neumann gave the only well-known first solution [23] for the logistic $x1 = 4(x - x^2)$. Therefore, two important directions are used in physics: quantum chaos [11] and thermodynamics [22].

This paper brings many improvements and corrections to the first preprints [3] [4] [5].
We apologize especially about some errors remaining in the previous versions of this paper.

***The subject of this paper is to describe the structures of invariant probabilistic distributions when we iterate a function indefinitely. They can exist, even they are masked by others cycles or distributions. The convergence to some point or some distribution depends on zones of domination, as we match with the topographic zones of attraction of certain cycles or fixed points. This problem is not completely studied here. The reader is warned that we put forward here a very static and asymptotic conception of the physics phenomena. No transitions, no dynamics, no orbits, except perhaps for the cycles, only a probability of presence somewhere after infinite iterations in a compact set.***

# Chapter 1: The resolving deviation

## I- Hypothesis, notations and definitions

### 1--H0- Hypothesis on the polynomial application $f$

Let $f$ be a polynomial application of a compact set $C$ of $R^d$ in $C$ with at least one fixed point $f(0) = 0$ well isolated in $C$. Diameter of $C$ is $D$. $f$ is locally invertible at $0$. $f = (f_1, f_2, ..., f_l, ..., f_d)$. Then, $f$ is $C^\infty$ and have a finite number $\nu$ of reciprocal images.

### 2- Collection of iterations

Let $C(f)$ be the collection of iterations $C(f) = (f, f^{(2)}, ..., f^{(k)}...)$ where $f^{(k)} = f \circ f^{(k-1)}$. We can index generically an $f^{(k)} \in C(f)$ as $f_\bullet$ with a spot « $_\bullet$ ». We check H0 for $f_\bullet$ provided to be invertible at $0_\bullet$ (see below). So, generally, all demonstrations for $f$ are valid for all $f_\bullet$.

### 3- Notations of the fixed points and eigenvalues of the linear part at theses points

Let $0 \in Fix(f) : f(0) = 0$. But, all $f^{(k)}$ can have $f^{(k)}(\alpha) = \alpha$ in $C$ with $f(\alpha) \neq \alpha$. These $\alpha_k$ are generally called cycles of $k$ order. Let $Fix(f^{(k)})$ be this set and $Fix(f_\bullet)$ the generic union for all $k$. Under H0, the real eigenvalues $\lambda$ ($\lambda_\bullet$) of the linear part of $f$ ($f_\bullet$) at $0$ ($0_\bullet$) are determining for the convergence to some invariants: fixed points, cycles or distributions. We don't study here complex or multiplicity of $\lambda$. These cases are not very difficult. When exist integers $r \in N^d$ such as $\lambda^r = 1$, there is resonance, and we present here only a trivial situation of this recurrent case. Let $\operatorname{Re} s^+ = \left\{ k \in N^d \left| |\lambda|^k > 1, |\lambda|^{k+1_l} > 1, l = 1, ...d \right. \right\}$ and $\operatorname{Re} s^- = \left\{ k \in N^d \left| |\lambda|^k < 1, |\lambda|^{k+1_l} < 1, l = 1, ...d \right. \right\}$

### 4- General notations

If $a, x, y \in R^d$, $xy$ is the scalar product; $x.y$ the vector with coordinates $x_l y_l$, $l = 1, .., d$; in **bold** $\boldsymbol{x} = x_1 ... x_l . x_d$; $d\boldsymbol{x} = dx_1 .. dx_l . dx_d$ and $\partial h / \partial \boldsymbol{a} = \partial^d h / \partial a_1 .. \partial a_l . \partial a_d$.

## II – Laplace-Steiltjes transform of the Perron Frobenius (PF) equation

### 1- Definitions : Perron Frobenius's operator (PF) and its Laplace-Steiltjes transformation (LPF)

Let $P$ be an arbitrary probability measure and $f$ an H0-application. For all Borelian set $B$, the transformation of $P$ by the Perron Frobenius's operator (PF) is $P \circ f^{-1}(B) = P_f(B)$. The measure is $f$-invariant if $P(B) = P_f(B)$ for all Borelian set $B$. As $f$ (resp. $f_\bullet$) apply a compact set $C$ in $C$, $P$ (resp. $P_\bullet$) exist $[15]$. All $f^{(k)} \in C(f)$ lets $P_f = P$ invariant.

### 2- Lemma 1

Under H0, for all $f$ and $p$, $\phi_f(y)$ is: $\phi_f(y) = \int_{R^d} e^{yf(x)} dP(x)$

By definition of $P_f(x) = P \circ f^{-1}(x)$, $\int h(x) dP_f(x) = \int h \circ f(x) dP(x)$ for all positive function $h$ of $L^1(P)$ (see $[15]$), hence for $e^{yx}$. Here, we note $y$ instead of $-y$ as traditionally to lighten the writing (especially for the approximations). There is no matter because we work in a compact set.

If $P$ has a Lebesgue-Steiltjes density $p$, then $p_f = \Sigma_\nu p \circ f_\nu^{-1} \left| f_\nu^{-1} {}' \right|$ is the transformation of $p$ by (PF) $[15]$. But this definition is not very useful, complicated and sometimes misleading. We prefer use the lemma 1. Nevertheless, if we suppose that we start the iteration process from a point $x_0$: either almost all points are in a neighbourhood of $0$ and the Lebesgue-Steiltjes density is the Dirac's distribution $\delta(x)$, or $x_n$ does not tends to $0$ with a distribution $p(x)$ vanishing at $0$.

Now, we take the Laplace-Steiltjes transform $L$ of a probability measure $dP(x)$ with support contained in the compact $C$. This transform is naturally used to capture the mass of the fixed point $0$:

$\phi(y) = \lim_{\varepsilon \to 0} \int_{-\varepsilon}^{+\infty} e^{yx} dP(x)$. Then, we have $\phi_f(y) = \lim_{\varepsilon \to 0} \int_{-\varepsilon}^{+\infty} e^{yx} dP_f(x)$ for all $f$ with a fixed point $0$ and all probability measure $P$ on $C$. Its existence is due to the compactness de $C$, which makes the set $\Omega = \{y \mid \lim_{\varepsilon \to 0} \int_{-\varepsilon}^{+\infty} e^{yx} dP_f(x) < \infty \mid \lim_{\varepsilon \to 0} \int_{-\varepsilon}^{+\infty} e^{yx} dP(x) < \infty \}$ not empty: $0 \in \Omega$ convex.

### 3- Remark: analytical representation of $\phi(y)$

Under H0, the support of the distribution is contained in the compact $C$ with diameter $D$. Hence, $\phi(y)$ can be represented by analytic series in a poly-disk with infinite radius of convergence. $\phi(y) = \Sigma_n b_n y^n$ converges for all probability measure $P$ in $C$ ($|b_n| < D^n / n!$). (More, we can take $y = k(n+1)$ provided $Dk e^{Dk} < 1/ed$ to keep the convergence of the series).

### 4- Definition: resolving equation L

With the Laplace's transform, the PF equation gives a resolving equation L:
$\theta_f(\phi) = \phi_f(y) - \phi(y) = 0$ and a collection of equations $\theta^\bullet_{f_\bullet}(\phi^\bullet) = \phi^\bullet_{f_\bullet}(y) - \phi^\bullet(y) = 0$ for $f_\bullet$.
More, if $P$ is invariant under $f$, $\phi_f(y) = \int e^{yf \circ f(x)} dP(x) = \phi_{f \circ f}(y)$. By recurrence: $\phi_{f_\bullet}(y) = \phi(y)$.

## III – Main Result

### 1- Notations of the sets of zeros of $e^n(y)$

Under H0, without resonance, let $\varphi(y) = \Sigma_n c_n y^n$ be a solution of the resolving equation. We note $\theta_f(\varphi) = \Sigma_n c_n e^n(y) = 0$ where $e^n(y) = (y^n - H_n(y))$ and $H_n(y) = \partial^n e^{yf(a)} / \partial a^n \mid a = 0$ at $a = 0$. Let $\varphi_n$ be the expansion of order $n$ de $\varphi$ such as $\theta_f(\varphi_n) = 0$.

*If $Z(e^n)$ is the set of zeros of $e^n(y)$, let $Z(E_n) = Z(e^n) \cap_{\forall l} Z(e^{n+1_k})$ $l = 1, \ldots d$, be the set of common zeros of these polynomials. $E_n$ is the ideal generated by these $d + 1$ polynomials.*

### 2- Theorem 1

*When $n \to \infty$:*

*1- Exists an integer $m$ such as $(\varphi_n(y) - 1)^m \in E_n$. With a probability 1, the distribution defined by $\varphi_n(y)$ is constant on the algebraic affine real manifolds defined by $E_n$.*

*2- If the manifolds are 0-dimensionnal and if $\varphi(s)$ is solution de the equation $\theta_f(\varphi) = 0$ for the density $q(a)$, we can write that the invariant density of Perron-Frobenius $p(a)$ is*

$$p(a) = \prod_{l=1}^{l=d} (-a_l) \partial^d q(a) / \partial a_1 .. \partial a_l .. \partial a_d$$

*3- These results are valid for $\forall f_\bullet \in C(f)$.*

The demonstration needs some steps. $\delta(x)$ is the Dirac's distribution at $0$;

## IV- Effects of a translation on analytical expression of $\phi_f(y)$

### 1- Notations of the effects of a translation $a$:

Observe the effects of a translation $a$: $x = u + a \neq 0$ changing the origin from $0$ to $-a$ in $C$. If $p(x)$ is a distribution density of a Lebesgue-Steiltjes measure for all $x \neq 0$, we have:

| | | |
|---|---|---|
| $f(x)$ | becomes | $f_a(u) = f(u+a) - a$, |
| $p(x)$ | becomes | $p(u) = p(x-a)$, |
| $p_f(x)$ | becomes | $p_f(x-a)$, |
| $\phi(y)$ | becomes | $\phi_a(y) = \phi(y,a) = \phi(y)e^{ya}$, |
| $\phi_f(y) = E(e^{yf(x)})$ | becomes | $\phi_f(y,a) = E(e^{yf(u+a)})$, |
| $\theta_f(\phi)$ | becomes | $\theta_f(\phi_a) = \phi_f(y,a) - \phi(y,a)$. |

The following proposition uses the notations: $H_n(y,a)e^{yf(a)} = \partial^n e^{yf(a)}/\partial a^n$ and $H_n(y) = H_n(y,0)$. All following series are uniformly convergent for all $y$ (see Constantine [6] with $|b_n| < D^n/n!$). We shall see uniform valuations later.

**2- Proposition 1**

- Let $\phi(y) = \Sigma_n b_n y^n$ be the Laplace-Steiltjes transform of an arbitrary $p$. Under H0, the iteration $f$ transforms all terms $y^n e^{ya} = \partial^n e^{ya}/\partial a^n$ of $\phi_a(y) = \phi(y)e^{ya}$ in $\partial^n e^{yf(a)}/\partial a^n$, and, for all $f$ with $f(0) = 0$ and $p$, we have:

$$\boxed{\phi_f(y,a) = \Sigma_n b_n \partial^n e^{yf(a)}/\partial a^n}$$

- Invariant density is solution of the resolving equation: $\theta_f(\phi) = \Sigma_n b_n(y^n - H_n(y)) = 0$. Moreover, for $\forall a \in C$ and $\forall y$, $\theta_f(\phi_a) = \Sigma_n b_n \partial^n(e^{ya} - e^{yf(a)})/\partial a^n = 0$. For all multi-integers $q$, we already have $\partial^q \theta_f(\phi_a)/\partial a^q = \partial^q E(e^{ya} - e^{yf(u+a)})/\partial a^q = \Sigma_n b_n \partial^{n+q}(e^{ya} - e^{yf(a)})/\partial a^{n+q} = 0$ and $\partial^q \theta_f(\phi) = \Sigma_n b_n(y^{n+q} - H_{n+q}(y)) = 0$.

- As $\phi(z) = L(p(x))$ is the Laplace's transform of the density $p$, we write $p(x-a) = L^{-1}(\phi(z)e^{za})$. As $\phi_f(y,a) = \int e^{yf(x)}dP(x-a) = \int e^{yf(x)}p(x-a)dx$, we have: $\phi_f(y,a) = \int e^{yf(x)}L^{-1}(\phi(z)e^{za})dx$
Here, $\phi(z)$ has the analytic expansion $\phi(z) = \Sigma_n b_n z^n$:
$\phi_f(y,a) = \int e^{yf(x)}L^{-1}(\Sigma_n b_n z^n e^{za})dx = \Sigma_n b_n(\int e^{yf(x)}L^{-1}(z^n e^{za})dx)$

$L^{-1}(z^n e^{za})$ is the product of $d$ one-dimensional integrals $L^{-1}(z^{n_l} e^{za_l})$ for $l = 1,...,d$. Using the Bromwich's formula, we recognize an inverse Fourier transform: $L^{-1}(z^{n_l} e^{za_l}) = F^{-1}(z^{n_l} e^{za_l}) = \delta^{n_l}(x_l - a_l)$ because $z^{n_l} e^{za_l}$ has no singularities (see Bromwich's formula on Wikipedia). $L^{-1}(z^n e^{za}) = \partial^n(\delta(x-a)/\partial x^n = \delta^n(x-a)$. Then:
$\phi_f(y,a) = \Sigma_n b_n(\int e^{yf(x)}\delta^n(x-a)dx) = \Sigma_n b_n \partial^n e^{yf(a)}/\partial a^n$

- Hence, the resolving equation is: $\theta_f(\phi_a) = \phi(y,a) - \phi_f(y,a) = \Sigma_n b_n \partial^n(e^{ya} - e^{yf(a)})/\partial a^n = 0$. Then, at $a = 0$, we have $\theta_f(\phi) = \Sigma_n b_n(y^n - H_n(y)) = 0$. In $\theta_f(\phi_a) = 0$, we observe the identity: $\phi(y) \equiv (\Sigma_n b_n H_n(y,a))e^{y(f(a)-a)}$ for $\forall a \in C$ and $\forall y$, and the last result is obtained by derivation. Denote that, in combinatorial theory, these polynomials $H_n(y,a)$ are known as Bell's polynomials. (Notice: $H_n(y)$ are polynomials in $y$ even if $f$ is only analytical).

**3- Remarks**

- We can localise the results on a $p$-dimensional manifold taking a volume on the manifold instead of the Lebesgue-Steiltjes measure on $R^d$.

- A diffeomorphism $g$ in $C$ of infinite class defined by $x = g(u)$ induce $f_g(u) = g^{-1} \circ f \circ g(u)$, hence $\phi_f(y)$ becomes $\phi_{g^{-1} \circ f \circ g}(y)$.

- If $a = \alpha_k^p \in Fix(f^{(k)})$ is a point of a cycle of order $k$ for $p = 1,...,k-1$ and $\alpha_k^{p+1} = f(\alpha_k^p)$, then $\phi(y) = (\Sigma_n b_n H_n(y,\alpha_k^p))e^{y(\alpha_k^{p+1} - \alpha_k^p)}$ at $\alpha_k^p$. So: $[\phi(y)]^k = \prod_{p=1}^{p=k}(\Sigma_n b_n H_n(y,\alpha_k^p))$.

**4- Corollary 1**

*Polynomials $H_n(y,a)$ have some recurrent matrix equations (see Constantine [6]):*

- $H_{n+1}(y,a) = \partial H_n(y,a)/\partial a + H_n(y,a)y\partial f(a)/\partial a$

- $\partial H_n(y,a) / \partial y = \Sigma_{k=0}^{k=n} C_n^k H_{n-k}(y,a) \partial^k f(a) / \partial a^k$.
- At $a = 0$, the leading term of $H_n(y)$ is $(\lambda y)^n$ and polynomial $\lambda^{-n} H_n(y)$ can be seen like « unitary » in $y$.

## V- Relation between solutions $\varphi(y)$ of (LPF) and $\phi(y)$ of (PF)

Let $\varphi(y)$ be the solution of (LPF) $\theta(\varphi) = 0$ and $\phi(y)$ the transform of the solution of (PF). For conveniences, we suppose the random variable $d$-dimensional.

**Proposition 2**

*Under H0, $\phi(y,a) = \partial^2 \varphi(y,a) / \partial y \partial a = \partial(y\varphi(y,a)) / \partial y$*

Let $\varphi(y)$ be a solution of $\theta_f(\varphi) = 0$. In the expansion of $\varphi(y)$, $b_n$ are fixed and we can construct $\varphi(y,a)$ and $\theta_f(\varphi_a) = 0$. But, according to proposition 1, if $\varphi(y,a)$ is solution of (PF), $\theta_f(\varphi_a) = 0$ will be an identity for all $y$ and $a$. All derivatives must vanish. In fine: $\theta(\partial^2 \varphi(y,a) / \partial y \partial a) = 0$. $\phi(y,a) = \partial^2 \varphi(y,a) / \partial y \partial a$. As $\varphi(y,a) = e^{ya} \int e^{ys} dP(s)$, we have the result.

It is not difficult to rewrite the result for a random vector $p$-dimensional $p \leq d$ on a manifold.

We can verify the arguments are true for every $f_\bullet$. In the following, no resonance gets very important.

## VI- Analysis of the solution $\varphi(y)$ of LPF

$\varphi(y)$ is now the Laplace-Steiltjes transform of a distribution with support in $C$ and solution of LPF. We study the expansion $\varphi_n(y)$ limited to the order $n$ and define the resolving deviation:

### 1- Definition 4: resolving deviation and the ideal $E_n$

*We call resolving deviation: $e^n(y,a) = \partial^n (e^{ya} - e^{yf(a)}) / \partial a^n$. At the fixed point $a = 0$, $e^n(y) = e^n(y,0) = y^n - H_n(y)$ is a polynomial of degree $n$.*

### 2- Proposition 3

*If $(1 \neq \lambda^n)$, under H0 and $\theta_f(\varphi_n) = 0$, we can construct $d$ polynomials $\varphi_n$ on $Z(E_n)$ and $\lim_{n \to \infty} Z(E_n) = (y | \theta_f(\varphi) = 0)$.*

For all analytic function $\psi(y)$, we can form the function $\theta_f(\psi) = \psi_f - \psi$ where $\psi_f(y)$ is obtained in substituting $y^n$ by $H_n(y)$ in the analytic expansion of $\psi$. We perturb a solution $\varphi$ of $\theta_f(\varphi) = 0$ by an arbitrary polynomial $\varepsilon_n(y) = \sum_{k=0}^{k=n} \varepsilon_k y^k$. The coefficient of the leading term of $\varepsilon_n(y)$ is: $\varepsilon_n \neq 0$.

The perturbed function $\psi^n(y) = \varphi(y) + \varepsilon_n(y)$ becomes under $f$: $\psi_f^n = \varphi_f + \sum_{k=0}^{k=n} \varepsilon_k H_k(y)$, yet convergent, and $\theta_f(\psi) = \sum_{k=0}^{k=n} \varepsilon_k e^k(y)$. As $\varepsilon_n \neq 0$ and $1 \neq \lambda^n$, the coefficient of the leading term of $\theta_f(\psi^n)$ is the one of the leading term of $e^n(y)$ if $y$ is non-null: $\varepsilon_n(1 - \lambda^n) \neq 0$. We cannot vanish this coefficient. We must seek a polynomial with a greater degree, and so on.

By induction, to vanish $\theta_f(\varphi)$, we must take $e^n(y)$ null when $n \to \infty$ for all non-null $y$.

Now, we take a very large fixed $n$ to construct $\varphi_n$ by this way. Suppose $y \in Z(e^n)$ for this $n \in N^d$.

Can we have $\varphi_n(y) = \Sigma_{p \leq n} c_p y^p$ as solution of $\theta_f(\varphi_n) = \sum_{p=0}^{p=n} c_p e^p(y) = 0$?

Response is yes: We find $d$ solutions, but under some other conditions.

Here, polynomial $e^n(y) = \sum_{k=0}^{k=n} \alpha_{k,n} y^k$ is known, but vanished if $y \in Z(e^n)$. $\theta_f(\varphi_n)$ has a smaller degree than $n$: we have $d$ candidates $e^{n-1_l}(y)$ for the new leading term. Choosing one $e^{n-1_l}(y)$, we must identify the coefficients $c_p$ to compute $\varphi_{1n}$. The term of degree $n-1_l$ is cancelled with $c_{n-1_l}$: $(1-\lambda^{n-1_l})c_{n-1_l} + c_n \alpha_{n-1_l,n} = 0$. Under the no-resonance hypothesis, notice that $1 - \lambda^{n-1_l} \neq 0$ for $l = 1,...,d$. But $c_n$ is arbitrary. Now, the polynomial $c_n(\alpha_{n-1_l} e^{n-1_l}(y) - (1-\lambda^{n-1_l})e^n(y))$ is known. By successive linear identifications, we shorten the degree of unknown coefficients of $\varphi_{1n}$. If we denote coefficients $c_k = c_{k,1}c_n$, we have $P_{n-1_l}(y) = c_n(\sum_{k \leq n} c_{k,1} e^k(y)) = 0$, then $\theta_f(\varphi_{1n}) = 0$. That is true for $\forall l$ and we have $d$ solutions $\varphi_{1n}$. We take one denoted $\varphi_n(y)$.

We examine $\varphi_n(y)$ with a translation $a$ of the origin: $\varphi_n(y)$ becomes $\varphi_n(y,a)$. As our previous computation for $\varphi(y)$, we have also for all $a$: $\theta_f(\varphi_n(y,a)) = \varphi_n(y)e^{ya} - \varphi_{fn}(y,a) = 0$. We expend $\theta_f$: $\theta_f(\varphi_n(y,a)) = \sum_{k \leq n} c_k e^k(y,a))$. Then, we can derive $\theta_f$ to obtain at $a = 0$: $\partial \theta_f(\varphi_n(y,a))/\partial a_k|_{a=0} = \sum_{m \leq n} c_m e^{m+1_k}(y) = 0$ for $k = 1,..,d$. A solution of this equation is obtained in taking the polynomial $e^{n+1_k}(y) = 0$. Hence, if $y \in Z(e^n) \cap_{\forall k} Z(e^{n+1_k})$, $\varphi_n(y)$ is solution of the both equations $\theta_f(\varphi_n) = 0$ and $\theta_f(\varphi_{n+1_k}) = 0$. Thence, the $d$ solutions $\varphi_{1n}$ are valid.

**3-Some inequalities**

*Let $\varphi$ be the solution of $\theta_f(\varphi) = 0$ and $\varphi_n(y)$ the expansion of order $n$. We have:*
$|\varphi_n(y) - \varphi_{fn}(y)| < 2d|yD|^{n+1} e^{|yD|}/(n+1)!$, $|\varphi(y) - \varphi_n(y,a)| < (e^{|ya|}-1)e^{|yD|} + |Dy|^{n+1} e^{|yD|}/(n+1)!$

The general term of this series is less than $|yD|^n/n!$. The remaining Lagrange-Taylor's rest of order $n$ of $e^{ya}$ is less than $|\varphi(y) - \varphi_n(y)| < d|yD|^{n+1} e^{|yD|}/(n+1)!$. As $|f(a)| < D$, we have also: $|\varphi_f(y) - \varphi_{fn}(y)| < d|yD|^{n+1} e^{|yD|}/(n+1)!$ for $\forall |a| < D$, and $|\varphi(y) - \varphi_n(y,a)| < (e^{|ya|}-1)e^{|yD|} + |Dy|^{n+1} e^{|yD|}/(n+1)!$.

Now, for fixed $n$ and $n - 1_l$, we study the solution $\varphi_n$, with $\varphi_n(0) = 1$ at $y = 0$.

**4- Lemma 2**

*For large $n$, $\varphi_{1n}$ converge to $\varphi$ when $n \to \infty$. Exist a number $m$ such as $(\varphi_n - 1)^m \in E_n$.*

We take $y \in \{Z(E_n) \cap y \leq k(n+1)\}$. As $\varphi_{1n}$ is a Laplace-Steiltjes transform of a measure in C, we have $|\varphi_{1n}(y) - \varphi_{1fn}(y)| < d|yD|^{n+1} e^{|yD|}/(n+1)!$, and then $\varphi_{1n}$ converge to $\varphi$ when $n \to \infty$. We can identify $\varphi_{1n}$ with $\varphi_n(y) = \varphi_n$. In such case, the coefficient of the leading term of $\varphi_n(y)$ and the one of $\varphi_{fn}(y)$ are close to 0 when $n \to \infty$. But, $\theta_f(\varphi_n) = e^n(y)$, and the leading term of $e^n(y)$ is $1 - \lambda^n$ and tend toward 1 or infinite. Either $\varphi_n - 1$ has same zeros than $Z(E_n)$ or $y = 0$. Imbedding the problem in the algebraically closed field of the complex numbers, the Hilbert's theorem of zeros (nullstellensatz) says $(\varphi_n - 1)^m \in E_n$.

Now, taking $y$ complex, we give a first characterisation of the distribution on $I_n$ as in Lucaks [17] or Wintner.

**5- Proposition 4**

*With a probability 1, all real point of $E_n$ belongs to the support of the distribution $L^{-1}(\varphi_n)$ and tends asymptotically to the support of $L^{-1}(\varphi)$. The distribution $p_n(a) = L^{-1}(\varphi_n)$, defined for $y \in Z(E_n) - \{0\}$, is constant on the real algebraic affine manifolfs (connected component) of $E_n$, null elsewhere.*

Let $y \in E_n \subset C^d$ and $\varphi_{1n}(y) = 1$. Consider the algebra quotient $A$ to define the affine manifold. If $y = v + iw$ has a real part $v \neq 0$, $v$ belongs to the real part of $A$. So we have $\varphi_n(y) = \int e^{ya} dP_n(a) = 1$. But, by definition of the Laplace transformation, we have the same result for the real algebraic affine manifold $\int e^{va} dP_n(a) = 1$. By difference: $\int e^{va}(\cos wa - 1) dP_n(a) = 0$. As $\cos wa - 1$ is negative continue, $P_n$ must be purely discontinue at $\cos wa = 1$ and $wa = 2k\pi; k \in Z$. We have $P_n(a) = \Sigma_{k \in Z} p_k Y(wa - 2k\pi), k \in Z$: a periodic lattice distribution. Here $Y$ is the Heaviside's function for the discrete distribution $p_k$ with sum $\Sigma_k p_k = 1$.

$\varphi$ verifies $\int e^{ya} dP(a) = \int e^{yf_\bullet(a)} dP(a)$ and for $\varphi_n$: $\int e^{ya} dP_n(a) = \int e^{yf_\bullet(a)} dP_n(a)$ whatever $y \in Z(E_n)$ when $n \to \infty$. But $\varphi_n = 1$ on $Z(E_n)$. So $wf_\bullet(a) \approx 2k\pi; k \in Z$. Consequently, $w(f_\bullet(a) - a) \approx 0$ for $\forall w \in \text{Im}(E_n)$. Then $a = 0$ and $p_0 = 1$ for $\forall v \in \text{Re}(E_n)$.

**6- Remarks**

- In resonant case, exist $n \in \zeta_R = \{r \in N^d \mid 1 = \lambda^r\}$ where the leading term of $H_n(y)$ is $y^n$. To construct $\varphi_n$, we can identify the $c_p$ for all $n \in \zeta_R$. But we shall see more interesting methods.

- If the measure $dP(a)$ is invariant, on have the identity $\int e^{ya} dP(a) \equiv \int e^{yf^{(p)}(a)} dP(a)$ for $\forall p \in N$. $\theta_p(y) = \phi_p(y) - \phi(y) = 0$. Then, the resolving deviation $f^{(p)}$ will asymptotically vanish.

**7- Lemma 3**

*If $n \to \infty$ with $n \in \text{Re } s^+$, the zeros of $e^n(y)$ are very close to theses of $H_n(y)$. Conversely, if $n \to \infty$ with $n \in \text{Re } s^-$ and $n \to \infty$, the solution is $y = 0$. We can say that the ideal $E_n - \{0\}$ is asymptotically generated by the zeros of polynomials $H_n(y)$ and $H_{n-1_l}(y) \; l = 1, 2, .., d$.*

If $n \in \text{Re } s^+$ and $n \to \infty$, we can write $e^n(y)$ as $e^n(y) = \lambda^n(\lambda^{-n} y^n - \lambda^{-n} H_n(y))$. But the leading term of $\lambda^{-n} H_n(y)$ is 1 and $\lambda^{-n} y^n \to 0$. By the continuity of the roots, we can say that the zeros of $e^n(y)$ are very close to theses of $H_n(y)$. Now, If $n \in \text{Re } s^-$ and $n \to \infty$, $e^n(y)$ is has its leading term of $y^n$, and $e^n(y)$ is very close to $y^n$. Solution is $y = 0$

**8- Normalisation of $E_n$**

The number of zeros can increase with $n$. For instance, if $k \in R$, suppose that we have $yf(ka)$ such $yf(ka) \to -\infty$ when $k \to \pm\infty$ for a fixed couple $(y, a)$, then, by Rolle's theorem $d^n e^{yf(ka)} / dk^n = \partial^n e^{yf(ka)} / \partial a^n \, a^n$ has $n$ zeros. We normalise these zeros (think to the Hermitian polynomials) to obtain a continuum. Let $y = n.s$ ($y_l = s_l n_l$): $E_n(y)$ becomes $\overline{E}_n = E_n(n.s)$.

**9- Lemma 4**

$\varphi_n(s) \leq 1$ *and* $Z(\overline{E}_n)$ *is bounded.*

- If $(\varphi_n(n.s) - 1) \in \overline{E}_n$, $\varphi_n(ks)$ must converge on $Z(\overline{E}_n)$. But logarithm $\varphi_n$ is convex. For every coordinate $s_l n_l$: $\log \varphi_n(s_l n_l) \leq \log \varphi_n(s_l(n_l - 1))/2 + \log \varphi_n(s_l(n_l + 1))/2$. Hence: $\log \varphi_n(s_l n_l) - \log \varphi_n(s_l(n_l - 1)) \leq \log \varphi_n(s_l(n_l + 1)) - \log \varphi_n(s_l n_l)$.

$\varphi_n(s_l n_l)) \geq \varphi_n(s_l) \varphi_n((n_l - 1) s_l) \geq (\varphi_n(s_l))^{n_l}$ and $\varphi_n(n.s) \geq \varphi_n(s)^n$ converge only if $\varphi_n(s) \leq 1$.

**VII- Random $d$-dimensional variables**

Let a $d$-dimensional density $q(s)$ such as $\varphi(y)$ is a solution of $\theta_f(\varphi) = 0$. Let $\phi(y)$ solution of the transform of PF for $f$ with density $p(s)$. Let $\partial C$ is the edge of the distribution.

**1- Proposition 5**

*If $q(s) = 0$ on $\partial C$, we have $p(s) = (-s) \partial q(s) / \partial s$*

From Proposition 2, $\phi(y, a) = \partial^2 \varphi(y, a) / \partial y \partial a$. Then, as $q(s) = 0$ on $\partial C$, we have the result.

**2- Hypothesis H1 on real zeros of $H_n(y) = \partial^n e^{yf(a)} / \partial a^n e^{-yf(a)} \big|_{a=0}$**

**H1-** Let $s_n \in S_n$ be such as $H_n(n.s_n) = 0$. The repartition limit of $s_n$ when $n \to \infty$ has a density $q(s)$ with respect to the $d$-dimensional Lebesgue-Steiltjes measure: in a volume $d\mathbf{y} = n^d d\mathbf{s}$ of $R^{+d}$ we have asymptotically $d\mathbf{y} = n^d q(s) ds$ zeros $s_n$.

Then, $n^d$ points $y_n$ verify $H_n(y_n) = 0$. Let $\varphi_n(y) = \frac{1}{n^d} \Sigma_n e^{yy_n}$

### 3- Proposition 6
*Under* H0 *and* H1, *if it is not masked, the Perron Frobénius's density is asymptotically*
$$p(a) = \prod_{\square=1}^{\square=d}(-a_\square)\partial^d q(a)/\partial a_1..\partial a_\square..\partial a_d$$

Let $t = y/n$. Under H0, H1 $\varphi_n(t.n) \approx \frac{1}{n^d}\Sigma_n e^{ts_n} n^d q(s_n) ds_n \to \varphi(t) = \int e^{ts} q(s) ds$.

## VIII- Applications

### 1- The logistic $f = \lambda(a - a^2)$ in dimension 1.
if $4 \geq \lambda \geq 1$, $f$ apply the compact set $C = [f(\lambda/4), \lambda/4]$ in C. Its resolving deviation is: $e_f^n(y)\big|_{a=0} = y^n - H_n(y)$ where $H_n(y) = H_n(\sqrt{\lambda y/2})(\sqrt{2\lambda y})^n$ with the Hermitian polynomials of order $n$ $H_n$. Hence, the distribution of the zeros of $H_n$, denoted as $t$, is asymptotically the semi-circular low of Wigner: $q(t)dt = k\sqrt{1-t^2} dt$. $q(1) = q(-1) = 0$. Thence $t = c\sqrt{y}/v$ with $c = \sqrt{\lambda/2}$ and $v = \sqrt{2n+1}$ has the semi-circular law. The distribution $p_1(s)$ of the variable $s = y/2n+1$ will verify: $t = r\sqrt{s}$ and: $p_1(s)ds = -(dq(t)/dt)dt = kr/\left(2s\sqrt{s(1-r^2 s)}\right) ds$. Then $p(s) = sp_1(s)$. The density of the invariant measure is a beta law $\beta(1/2, 1/2)$: $p(s) = k'/\sqrt{s(1-r^2 s)}$, obtained by Kawamura with the Cunts-Krieger's algebra, which generalise the well-known result of Ulam and Von Neumann. But this distribution will be masked if exist attractive cycles according with the values of $\lambda$.

### 2- First approach of the Julia's iteration
We prevent the readers that this example doesn't pretend characterise the Julia's sets. Many important works have been made. Here, we don't worry with cycles or resonances, or relations of domination. Let $f(z) = z^2 + 1/4 - \alpha^2$ be the iteration de Julia in C [14]. It has 2 fixed points: $a_0 = 1/2 + \alpha$ and $a'_0 = 1/2 - \alpha$. We denote here with « **bold** » characters the vectors of $R^2$ for the complex numbers $\mathbf{a} = (a,b)$, $\boldsymbol{\alpha} = (\alpha, \beta)$ and $\mathbf{a'_0} = (a'_0, b'_0) = 1/2 - \boldsymbol{\alpha}$.

Consider the function $f(z)$ at the fixed point $\mathbf{a_0}$; we can write with $z = a + a_0$: $\mathbf{f} = (a^2 - b^2 + 2a_0 a - 2b_0 b, 2ab + 2a_0 b + 2b_0 a)$. We can suppose that it applies a compact C in C by convergent reciprocal iterations. We have the function: $\mathbf{y}f(\mathbf{a}) = 2\,\mathbf{y}L(\mathbf{a}) + \mathbf{y}Q(\mathbf{a})$.

Let $\mu_0 = \sqrt{a_0^2 + b_0^2}$, $\mu = \sqrt{x^2 + y^2}$, $\gamma = \sqrt{(\mu+x)/2\mu}$ and $\delta = \sqrt{(\mu-x)/2\mu}$. We denote the orthogonal matrices:

$$\mu_0 \mathbf{O}_0 = \begin{pmatrix} a_0 & -b_0 \\ b_0 & a_0 \end{pmatrix} \qquad \mu \mathbf{O}_y = \begin{pmatrix} x & y \\ y & -x \end{pmatrix} \qquad \mathbf{O} = \begin{pmatrix} \gamma & -\delta \\ \delta & \gamma \end{pmatrix}$$

- The linear part is $\mathbf{y}L(\mathbf{a}) = x(a_0 a - b_0 b) + y(a_0 b + b_0 a) = \mu_0 \mathbf{y} \mathbf{O}_0 \mathbf{a}$ where $\mathbf{O}_0$ is orthogonal. We suppose that the eigenvalues of the linear part $a_0$ and its conjugated $\overline{a_0}$ have a module greater than 1.

- The quadratic part is orthogonal $\mathbf{y}Q(\mathbf{a}) = xa^2 - xb^2 + 2yab = \mu\,\mathbf{a}\mathbf{O}_y \mathbf{a}$. $\mathbf{O}_y$ has a characteristic equation $x^2 + y^2 - \mu^2 = 0$ with 2 eigenvalues, one positive $\mu$, and the other negative $-\mu$.

If $^tO$ design the transposed matrix of $O$, denoting $y = Os$ and $a = Ou$ with $u = (u, v)$ and $s = (r, s)$, then $ya = su$. $yf(a) = 2\mu_0 \, s \, {}^tOO_0Ou + \mu u^2 - \mu v^2 = 2r(a_0 u - b_0 v) + 2s(a_0 v + b_0 u) + \mu u^2 - \mu v^2$

Hence: $yf(a) = (2ra_0 + 2tb_0)u + \mu u^2 + (2ta_0 - 2rb_0)v - \mu v^2$. In the orthogonal eigenspaces, the application splits into two independent iterations: $(2ra_0 + 2sb_0)u + \mu u^2$ and $(2ta_0 - 2rb_0)v - \mu v^2$.

As $O$ is orthogonal $\mu = \sqrt{r^2 + s^2}$.

The resolving deviation gives two Hermitian polynomials: $H_n((2ra_0 + 2sb_0)/i\sqrt{2\mu})(i\sqrt{2\mu})^n$ and $H_n((2sa_0 - 2rb_0)/\sqrt{2\mu})(\sqrt{2\mu})^n$. The first is always positive. The zeros of the second follow asymptotically a semi-circular Law. Thence, $(2s_1 a_0 - 2r_1 b_0)/\sqrt{2\mu}$ (with $r_1 = r/\sqrt{2n+1}$ and $s_1 = s/\sqrt{2n+1}$) follows a semi-circular Law. And $s_1\sqrt{2n+1} = s$ gives the distribution.

We change coordinates for $a'_0$ without modify the distribution. But, they may be masked (dominated) locally by every cycles or distributions of $f_\bullet$.

# Chapter 2: Approximation of the zeros

## I- Notations, rappels and hypotheses
### 1-Leading idea

We search the asymptotic distribution of the real zeros of $H_n(y)$ when $y \to \infty$ with $n$. As Plancherel and Rotach[20] did for the Hermitian polynomials, we represent $H_n(y)$ in the field of the complex numbers with the Cauchy's $d$-dimensional integral : $H_{n-1}(y) = \partial^{n-1} e^{yf(a)} / \partial a^{n-1}\big|_{a=0} = c \oint_\Gamma \frac{e^{yf(a)} da}{a^n}$ where $\Gamma$ is a closed poly-disk around the fixed point 0 of $f$, $a \in \mathbb{C}^d$. Note $n - 1 = (n_1 - 1, n_2 - 1, \ldots n_d - 1)$. $c$ is a finite non-null constant.

The steepest descent's method, due to Riemann-Debbie, gives us an approximation of this integral when $y \to \infty$. If $y$ is a zero of $H_{n-1}(y)$, we will hope that the approximation will give this zero. The method needs many assumptions. If $yf(a)$ is polynomial in $a$, we can see Pham [19] or Delabaere [7] for a good formulation, but only if the Hessian of $yf(a)$ is definite.

### Definition 1
*We call Plancherel-Rotach's function (PRF): $\gamma(a) = yf(a) - n\log a$ where $\log a$ is the main determination of the complex logarithm. Then $H_{n-1}(y) = c \oint_\Gamma e^{\gamma(a)} da$*

### 2-Basic notations

$y \in \mathbb{R}^d$ tends toward the infinite with $n \in N^d$ and we write $y = n.s$ where $y_i = n_i s_i$. All the $n_i \to \infty$. We shall see soon that all the $n_i$ must be equal. Now we can write $\gamma(a) = n.(sf(a) - \log a)$

### 3-Recalls on the steepest descent's method and assumptions

The approximation of $H_{n-1}(y)$ by the steepest descent's method is the finite sum of contributions $\Sigma_c Q_c(a)$ of all the critical points of $\gamma(a)$. Among these contributions, some of them are exponentially negligible. So, we keep only that one of which the real part is the most important. The approximation can be written (see [7]):

$$H_{n-1}(ns) = Q_n(a) = c' e^{\gamma(a)}(1 + g(a,n,s)/n')$$

Where the function $g$ is bounded when $n'$ depends on $n$ and $n' \to \infty$ with $n$. $c'$ is never null.
A point $a$ is critical point of $\gamma(a)$ if we have $\partial \gamma(a)/\partial a = 0$. That can be written $a_i.n.s \partial f(a)/\partial a_i - n_i = 0$. Then $a \partial \gamma(a)/\partial a$ is polynomial with the same degree than $f$.

### 4- Hypotheses: *general position*

*The critical point $a$ giving the maximum of $e^{\gamma(a)}$ is assumed be in **general position** [7]:*

### Definition 2
**H1**- *For all fixed $s$, it is isolated, unique, (except its complex conjugate) with finite distance. It is called Morse's point, if it is not degenerated.*
**H2**- *A sufficient condition to obtain this maximum is that the Hessian (which is Hermitian) of $\gamma(a)$ is definite negative at $a$. Critical points not degenerated are isolated.*
Unfortunately, in a lot of problems, the Hessian is not definite negative (see V).
**H3**- *$n.sf$ tends toward $-\infty$ when $n \to \infty$ in order to define the Laplace's integral.*
But, if the iteration is in a compact set, this hypothesis will be redundant.

## II- Asymptotic equivalence relation of the real zeros of $H_{n-1}(y)$ when $n \to \infty$

We first identify the equation of the zeros of the approximation, then we study under what conditions the distribution of the approximation is asymptotically very near of the zeros of $H_n(y)$. The variety of the situations is so rich that we examine here only the most common ones.

As we are only interesting by the real zeros of $H_{n-1}(y)$, we can neglect in the representation $Q(a)$ all the factors $\neq 0$ : that may be functions, constants, etc. Notice that $g/n' \to 0$ when $n \to \infty$ and $c' \neq 0$ in $Q(a) = c' e^{\gamma(a)}(1 + g(a, n. s, 1/n')/n')$. Then $e^{\gamma(a)}$ is the one and only one factor making $Q(a)$ null. We write this equivalence $H_{n-1}(y) \underset{0}{=} e^{\gamma(a)}$. Then $n'$ has no importance when $n \to \infty$.

**1-Lemma 1**

*Almost everywhere (when the Hessian is not degenerated):* $H_{n-1}(y) \underset{0}{=} J_n(y) = \sin(\operatorname{Im} \gamma(a))$

- First, if $H_{n-1}(y) = 0$, all contributions $e^{\gamma(a)}$ must be null at the critical point $a = a(y)$. As the polynomial $a \partial \gamma / \partial a = 0$ has real coefficients, if $a$ is a complex solution, the conjugate $\overline{a}$ is solution. But the contour gives a negative contribution for $\overline{a}$, and the sum is $e^{\gamma(a)} - e^{\overline{\gamma(a)}} \underset{0}{=} J_n(y)$.

- Conversely, if $J_n(a) = 0$, we can imagine that another contribution $Q_1(a)$ takes the place of $Q$. Hence, it can exist $y$ and $Q_1$ such as $|H_{n-1}(y)| = |Q_1(a)| = \delta > 0$ when $Q(a) = 0$. But, it is impossible because the two contributions are continuous:

Let $a'$ be a point in a small neighbourhood of $a$ where $Q$ is the dominating contribution and suppose that exists $Q_1$ dominated. At $a'$, $Q_1(a') < Q(a') < \varepsilon/3$ because $Q$ is dominating. Then $H_{n-1}(y') = Q(a')$. As $(a, a')$ determines a unique couple $(y, y')$ : $|H_{n-1}(y)| < |H_{n-1}(y) - H_{n-1}(y')| + |H_{n-1}(y')|$. In this inequality, we can substitute $H_{n-1}$ by $Q$ or $Q_1$ each time they are dominating: $|H_{n-1}(y)| < |Q_1(a) - Q(a')| + |Q(a')| < |Q_1(a) - Q(a')| + \varepsilon/3$, but $|Q_1(a) - Q(a')| < |Q_1(a) - Q_1(a')| + |Q_1(a') - Q(a')|$. The continuity of $Q_1$ induces the contradiction.

**2- Come back to the ideal $E_n$**

Recall that $Z(E_n)$ be the set of common zeros of the polynomials $e^n(y) = 0$ and $e^{n+1_l}(y) = 0$, $l = 1, \ldots d$. $E_n$ is ideal of these polynomials. We have seen that, if $n \to \infty$ with $n \in \operatorname{Re} s^+$, the zeros of $e^n(y)$ are very close to those of $H_n(y)$. But, we can study $H_n(y)$ even though $n \notin \operatorname{Re} s^+$. Let $Z(H_n)$ the set of common zeros of $H_n(y)$ instead of $e^n(y)$.

Let $Z(J_n)$ be the set of the common zeros of $J_n(y)$ for $n-1$ and $n-1+1_l$, $l = 1, 2, \ldots, d$.

**III – Main theorem**

Let $p$ be the fixed number of complex coordinates of a critical point $a$ of PRF.

**1- Theorem 2**

*If the critical point $a \in Z(J_n)$ is in general position, then:*

**1 -** $\gamma(a)$ *depends only on one common index of derivation* $n_1 = n$. *We can write* $y_1 = ns_1$ *and the PRF is* $n\gamma(a) = n(sf(a) - \log a)$ *with* $a = \Pi_1 a_1$. *The coordinates are defined by the system :*

$a_l s \partial f(a) / \partial a_l = 1$  $l = p+1, 2, \ldots, d$. *As $f$ is polynomial, solutions are algebraic affine manifolds. Moreover, the resolving deviation $e^n(y)$ must be computed only for common derivation $n$. Thus, let $\Delta = \lambda = |f|$ be the determinant of the linear part of $f$. We have*

*If $|\Delta| < 1$ the iteration converges to 0;*

*If $|\Delta| > 1$ the iteration, if it is not masked, can converge to a distribution defined by $Z(J_n)$*

**2 -** *The imaginary parts of the complex coordinates of $a$ satisfy a uniform code of the real zeros of $H_n(y)$, which are algebraic affine manifolds for all complex coordinates $l = 1, 2, \ldots, p$. These $d - p$ dimensional many folds are indexed by a code $k \in [1, n-1]^p$ and verify:*

$s_l \, \text{Im} \, f_l(a) - \vartheta_l = k_l \pi / n$ with the code $k_l \in (1,...,n-1)$. Between two zeros (manifolds $M_f$) at time $n-1$ there is one zero at time $n$: then we have a Rolle's foliation.

When $n_l \to \infty$

if $\kappa_l = \lim_{n \to \infty} k_l / n$ for $l = 1,..,p$, then $\kappa$ is asymptotically a random uniform vector on the unit cube $[0,1]^p$. Then, The $d-p$ dimensional many folds are indexed by a random vector of $[0,1]^p - \{0\}$ and define the random equation $\boxed{\pi \kappa_l = s_l \, \text{Im} \, f_l(a) - \vartheta_l}$  $l = 1,2,..,p$

**Remark**: It seems possible to extend the results to the rational fractions and analytical functions.

## 2- Proof: Analysis of the zeros of $J_n(y)$

Let $p \neq 0$ be the fixed number of complex coordinates of the critical point.
We must find the zeros of $\sin(\text{Im} \, \gamma(a))$. Hence $n \cdot (\text{Im} \, sf(a)) - \vartheta) = k \pi$.

### -Proof of point 1

Suppose we found a complex critical point $a(y)$ such as $J_n(y) = 0$. As $\sin k \, \text{Im} \, \gamma(a) = 0$ for all integer $k \in N$, $a$ is invariant if we multiply $y$ by any integer $k$. Hence, $k \gamma(a)$ gives the same critical point $a$ than $k' \gamma(a)$. We can write $k \gamma(a) = k' \Sigma_l n_l (s_l f_l(a) - \log a_l)$

But, if $y \in Z(J_n)$, all change of $n_l - 1$ in $n_l$ gives the same point $y$ for all $l$:
then $a_l y \partial f(a) / \partial a_l - n_l = 0$ becomes $(a_l + \Delta a_l) y \partial f(a + \Delta a) / \partial a_l - n_l + 1 = 0$. Dividing $y$ by a common $n = \min(n_l)$, we obtain by difference $\left| [\partial(a_l(y/n) \partial f(a) / \partial a_l) / \partial a] \right| |\Delta a| < 1/n$ and $|\Delta a| \leq c/n$. Then, the integrand is quite constant:

$k \gamma(a) = k' \Sigma_l (n_l + 1)(s_l f_l(a) - \log a_l) = k' \gamma(a + \Delta a) + k' \Sigma_l (s_l f_l(a + \Delta a) - \log a_l + \Delta a_l)$

$\left| (k - k') \gamma(a) - k' \Sigma_l (s_l f_l(a) - \log a_l) \right| < Cte / n$

Then, with $k - k' = 1$ all the $n_l$ must be equal: $n_l = n$ for $l = 1, 2,..,d$.

### -Lemma 2

If $k_{cl} = k_l / n$, the zeros of $J_n(a)$ must verify the condition $\Sigma_l n_l (s_l \, \text{Im} \, f_l(a) - \vartheta_l - k_{c_l} \pi) = 0$

$s \, \text{Im} \, f(a)$ contains only sinus, the arcs of which are linear combinations of multiple integers of $\vartheta$, then $s \, \text{Im} \, f(a)$ is periodic: $\text{Im} \, \gamma(a(\vartheta + 2k'\pi)) = \text{Im} \, \gamma(a(\vartheta)) + 2k'\pi$. Hence, all $N \, \text{Im} \, \gamma(a(2k'\pi)) - k\pi = 0$ are solutions. Writing $k\pi = \Sigma_l n k_{c_l} \pi$ where $k_{cl} = k_l / n_l$ and $y_{cl} = s_l n$, the solutions can be written as $n \Sigma_l (s_l \, \text{Im} \, f_l(a) - \vartheta - k_{c_l} \pi) = 0$ where $p$ is the number of complex coordinates. Let $k_l$ runs over $1, 2,..,n-1$, then the solution compounds $n-1$ distinct points for each complex coordinates.

In particular, $\vartheta = 0$ such as $\text{Im} \, \gamma(a(0)) = 0$ gives the fixed point and a non-null solution $\Sigma_l (s_l \, \text{Im} \, f_l(a) - k_l \pi) = 0$. But we don't know if there is other $0 \neq \vartheta \in (-\pi, +\pi)^p$ and under what conditions such a solution exists. As $a \in Z(J_n)$, we examine now a complex coordinate $a_l$, when $n_l - 1$ becomes $n_l$, the others $k \neq l$ remaining constant.

### - Proof of point 2

If $y \in Z(J_n)$, then $y$ remains invariant when $n_l - 1$ becomes $n_l$. The equation of the complex $a_l$, $a_l y \partial f(a) / \partial a_l = n_l$ becomes $a_l y \partial f(a) / \partial a_l - n_l + 1 = 0$, the other equations don't move. The critical point $a$ becomes $a + \Delta a$. If we normalise the equations by the common $n$, we have again:
$[\partial(a_l y \partial f(a) / \partial a_l) / \partial a] \Delta a = \Delta$ and $|\Delta a| < c \Delta$. Now, the equation
$\Sigma_j n_j (s_j \, \text{Im} \, f_j(a) - \vartheta_j - k_{c_j} \pi) = 0$ becomes: $\Sigma_j n'_j (s_j \, \text{Im} \, f_j(a') - \vartheta'_j - k'_{c_j} \pi) = 0$
where $n'_j = n_j$ if $j \neq l$ and $n'_l = n_l + 1$

$k'_{cj} = k_{cj}$ if $j \neq 1$ and $k'_{c1} = n_1 k_{c1}/(n_1+1)$ for the range of $k_{c1}$ + a new point : $n_1/(n_1+1)$

By difference we obtain; $\partial S(a+\mu\Delta a)/\partial a \Delta a + s_1 \operatorname{Im} f_1(a+\Delta a) - k'_{c1}\pi - \delta k'_{c_1}\pi = 0$

where $n_1 \delta k_{c_1}$ is the point when the range change from $(1, n_1 - 1)$ to $(1, n_1)$

Examine the quantity $\partial S(a+\mu\Delta a)/\partial a \Delta a$. As $\partial S(a)/\partial a = 0$ at the critical point $a$, then $|\partial S(a+\mu\Delta a)/\partial a| < |\partial^2 S(a+\mu\Delta a)/\partial a^2||\mu\Delta a|$. We obtain:

$|\partial \operatorname{Im} S(a+\mu\Delta a)/\partial a \Delta a| < c'|\Delta a| < c'c/n$

For all $n_1 k_{c_1}$ running over $(1, n_1)$, $s_1 \operatorname{Im} f_1(a) - \theta_1 - k_{c1}\pi \to 0$ when $n \to \infty$.

This is true for all complex coordinates and, when $n$ is large, we have the result.

We obtain a uniform code by $k_c$ of coordinates of real zeros of $e^{C(a)}\sin S(a)$. Note that a $n+1$-th zero is between two $n$-th zeros. If $0$ is the unique solution of $\operatorname{Im}(sf(a)) - \vartheta = 0$, the order $n+1$ manifolds can't intersect order $n$ manifolds. Then we have a Rolle's foliation

Suppose now that $n \to \infty$

We have already $s_\square = \lim_{n\to\infty} y_\square/n_\square, \kappa_\square = \lim_{n\to\infty} k_\square$ for $\square = 1, .., p$ ;

Let $\kappa$ a random vector with uniform repartition on the unit cube equation $s_\square \operatorname{Im} f_\square(a) - \vartheta_\square = k_\square \pi/n_\square$) tend to $s_\square \operatorname{Im} f_\square(a) - \vartheta_\square = \kappa_\square \pi$

### 3- Examples
### - One-dimensional case

$a_c$ depends only on $s = y/n : a_c \partial f(a_c)/\partial a - 1/s = 0$ and the equation of zeros is $(s \operatorname{Im} f(a_c) - \vartheta)/\pi = \kappa$. then $q(s)ds = \operatorname{Prob}\{1 \text{ zéro} \in (s, s+ds)\} = \operatorname{Im} f(a_c) ds/\pi$

### - Case of a cycle

Let a cycle : $f^{(p)}(0) = 0$ with $f^{(k)}(0) = \alpha_k \neq 0$ pour $k = 1, 2, .., p-1$. A point $\alpha_k$ of this cycle is fixed point for the application $f^{(p)}$ as $f^{(np)}(\alpha_k) = \alpha_k$. If exists a critical point $a_{ck}$ of $\gamma(a) = \zeta s f^{(p)}(a) - \zeta \log a$ we can have a probability distribution. But this distribution will move at each iteration and be cyclical as in a Markov process.

### - Hermitian case : $f = \lambda a - a^2/2 : c\mathrm{H}_{n-1}(\lambda) = c\partial^{n-1} e^{(\lambda a - a^2/2)}/\partial a^{n-1}|_{a=0} = \oint_\Gamma \dfrac{e^{(\lambda a - a^2/2)} da}{a^n}$

The PR function [11] gives for $t = \lambda/\sqrt{4n} = \cos\vartheta$ : $\mathrm{H}_{n-1}(\lambda) \underset{0}{=} \sin[\vartheta - \sin(2\vartheta)/2]$ when $\lambda \to \infty$ under the condition $\lambda^2 - 4n \leq 0$. The density of zeros of $t$ is the semi-circular low : $p(t)dt = (1 - \cos(2\vartheta))d\vartheta/\pi = (2/\pi)\sqrt{1-t^2} dt$.

Demonstration that seems more explicative than the theorem of Sturm.

### - r-Hermitian case: $f = \lambda a - a^r/r$.

We write the PR function $\gamma(z) = \lambda z - z^r/r - t \log z$ and the critical point is solution of the Lambert's [9] equation of $z\gamma'(z) = \lambda z - z^r - t = 0$.

If $r = 3$ and if we use the Fourier's transform, the characteristic function is the Airy's function.

### - Quadratic case in $\mathrm{R}^d$ : $f(a) = \lambda a + Q(a) : Q(a)$ is quadratic. We can diagonalize $yQ$ with transformation orthogonal $a = O_y u$. In this basis, the volume is invariant, then $n \log a = n \log u$. The function $\gamma(a)$ is $\gamma(u) = y\lambda O_y u_1 + \Sigma(\mu_1 u_1^2 - n \log u_1)$, the eigenvalues $\mu_1$ are solutions of $|yQ(a) - \mu I| = 0$. The semi circular low is generated by each coordinate having a $\mu_1 < 0$. Moreover, the distributions are independent in orthogonal eigenspaces. The equations $\mu_1 = 0$ give zones where the number $p$ of complex coordinates remains constant.

## IV- Important questions to be solved
### 1- Complexity or reality of the coordinates of the critical point

The question is to know how vary the number $p$ of complex coordinates of the critical point $a$ with $s$. We know some results in the 1-dimensional case. Criterions are condition as $0 < s \leq c"z$. In the multidimensional case, we obtain criterions with the Morse's theory [18], which can explain the Stokes's phenomena.

### 2- Relations of domination between dominating and masked distributions

- We have now a tool to explicit the relations of domination between the distributions imputed to each point of $Fix(f)$. Consider without loss of generality two fixed points: $\alpha$ for $f$ and $\alpha_m$ for $f^{(m)}$:

For fixed $n$, a point $s$ of a manifold $M_f$ solution at $\alpha$ verifies the $d + p$ following equations:

$\partial \gamma_f(a) / \partial a_l = 0$, $l = 1, 2, .., d$

$\operatorname{Im} \gamma_{fl}(a) = k_l \pi / n$, $l = 1, 2, .., p$

The first $d$ equations are generally a bijection between $a$ and $s$, The $n^p$ (for each $k_l$) others imply $n^p$ manifolds $d - p$ dimensional $M_f$.

We can do the same considerations to compute a point $s_m$ on a manifold $M_{f^{(m)}}$ solution at $\alpha_m$, but with $p_m$ instead of $p$.

Consider now the intersection of all these manifolds $M = M_f \cap M_{f^{(m)}}$, supposing we are in the Rolle's foliation case. We can write now $\gamma_f(s)$ instead of $\gamma_f(a)$ on $M$.

At a point $s$, one of 2 functions, $\exp(\gamma_f(s))$ or $\exp(\gamma_{f^{(k)}}(s))$, must exponentially dominate the other. Then, the dominated iteration is masked. A very important case for classify the region with constant domination is the surfaces defined by $\gamma_f(s) = \gamma_{f^{(k)}}(s)$ We can classify in this category the Stokes phenomena or all random manifolds like the Julia's iteration.

We shall apply this method of domination to EDP or EDO.

## V- Linear degeneracy of the Hessian

The Hessian of $y\partial^2 f(a) / \partial a \partial \bar{a}$ must be definite at the critical point for apply the previous results. But $\left| y\partial^2 f(a) / \partial a \partial \bar{a} \right|$ will degenerate in two distinguishable situations: this comes for particular values (from time to time) of $y$ or structurally, when the rank of $y\partial^2 f(a) / \partial a \partial \bar{a}$ is less than $d$.

### 1- The "from time to time" degeneracy

The determinant is a homogeneous polynomial in $y$ $\left| y\partial^2 f(a) / \partial a^2 \right|$ with a degree $d$. This polynomial can vanish when one or more eigenvalues vanish. This equation cut the space in regions where the eigenvalues of the Hessian keep a constant sign. For instance, in dimension 2 without loss of generality, if $f$ is quadratic, with two quadrics: $Q = Q(\alpha, \beta, \gamma) = \alpha a^2 + 2\beta ab + \gamma b^2$ and $Q' = Q(\alpha', \beta', \gamma')$. The Hessian of $yf$ where $y = (x, y)$: $xQ + yQ'$ has the determinant $(\alpha\gamma - \beta^2)x^2 + (\alpha'\gamma' - \beta'^2)y^2 + (\alpha\gamma + \alpha'\gamma' - 2\beta\beta')xy$. Its cancellation gives two straight lines of from-time-to-time degeneracy.

### 2- the structural degeneracy

More difficult is the structural degeneracy when, at the critical point, $\gamma(a)$ has a bloc of coordinates which depends on the others. We use here a **bold** notation of the vectors of $R^d$ with two spaces $R^p$ and $R^{d-p}$: $\boldsymbol{a} = (a, b)$ with $a \in R^p$ and $b \in R^{d-p}$. The structural degeneracy appears each time when $b$ depend functionally on $a$, without $a$ depends on $b$ in the equations $\partial \gamma(\boldsymbol{a}) / \partial \boldsymbol{a} = 0$. The Hessian is now functionally degenerated. Then the Hessian is function only of $a$ of $R^p$. The method of the steepest descent must be modified. Unfortunately, we don't know a good book about that.

Unlike the literature on the chaos, which imputes these phenomena the non-linearity, linearity will induce many complicated situations, or exhibit some auto-similar structures and many effects masked by the probability's calculus. The most elementary case studied here is that where $f$ is partly linear. This case gives a first explication of the " esthetical" chaos $[10]$.

### 3-Definition 3
The *iteration $f$ is said partly linear if $f$ is linear in $b$ and non linear in $a$*.

**Notations**

We can write that: $f(a) = f(a,b) = h(a)b + g(a)$ where $g$ is a polynomial application of $R^p$ in $R^d$ and $h$ is a polynomial matrix $(d, d - p)$ in $a$ each element of which apply $R^p$ in $R^d$. Here polynomial means non linear. At the fixed point $0 : g(0) = 0$. The linear part of $f$ at $0$ is supposed diagonalized real: $h(0) = (0, \lambda')$ and $\partial g(0)/\partial a = (\lambda, 0)$. We can write $f(a) = (A(a)b + B(a), C(a)b + D(a))$ where the $A(a)$, $C(a)$, have at least degree 0 and $B(a), D(a)$ degree 1. Note $a = (a,b)$ and $n = (n,n')$, $y = (x,y)$, $ay = ax + by$ and $yf(a) = yh(a)b + yg(a)$. The function of Plancherel-Rotach in $a$ is:

$n\gamma(a) = yh(a)b + yg(a) - n\log a - n'\log b$.

The convergence is deeply modified: the density for the d- dimensional Lebesgue's measure is null.

### 4- Proposition 3

*If $sh(a)$ has no real zero, the zeros of $H_{n-1}(y)$ tend to a family of algebraic manifolds of dimension p depending on an random uniform vector of $R^p$ for $\square = 1,2,...,p$ :*

$$\pi k_\square = s_\square \mathrm{Im}\, g_\square(a) + 1 Arg(sh_\square(a)) - \vartheta_\square$$

*When, $n \to \infty$ the code $\kappa_1 / n$ tends asymptotically to a distribution uniform on $(0,1)^p$ and gives the distribution of $s$ by reciprocal image:*

*The main manifold contains par the fixed point, the others are obtained par iteration of the critical point.*

*An important case is that where $h(a)$ is constant:* $\pi k_\square = s_\square \mathrm{Im}\, g_\square(a) - \vartheta_\square$

The resolving deviation can be written: $e_f^{n-1,n'}(a,y)_{a=0,b=0} = x^{n-1}y^{n'} - H_{n-1,n'}(y)$

with $H_{n-1,n'}(y) = \partial^{n-1}((yh(a))^{n'} e^{yf(a)})/\partial a^{n-1}\big|_{a=0,b=0}$ where the term $(yh(a))^{n'}$ is obtained by $n'$ derivations in $b$. Le polynomial $yh(a)$ is holomorphic at $a$. We can construct log$yh(a)$ on a locally connected open space, but we have many problems with the singularities of $yh(a)$, (number finite of zeros of $yh(a)$ of null measure) and the monodromy.

We represent $H_{n-1,n'}(y)$ by the integral of contour $\Gamma_a$ in $R^p$ :

$cH_{n-1,n'}(y) = \oint_{\Gamma_a} (yh(a))^{n'} e^{yf(a)}/(a^n)da = \oint_{\Gamma_a} e^{\gamma(a)}da$ where $b = 0$ in $f(a)$

$\gamma(a) = n'\log(yh(a)) + yg(a) - n\log a$.

Taking $n' = n = n - 1$, we note $y = ns$ and $1 \log a = \sum_{l=1}^{l=p} \log a_l$ and $1\log sh(a) = \sum_{\square=p+1}^{\square=d} \log sh_\square(a)$.

the PRF : $n\gamma(a) = n(1\log(sh(a)) + sg(a) - 1\log a)$

If the Hessian of $\gamma(a)$ is definite negative and if $sh(a)$ don't have real zero, the critical points $a$ verify the polynomial equations $\partial\gamma(a)/\partial a = 0$.

We take $\gamma(a)$ instead of the PRF $\gamma(a) = yf(a) - n\log a$. All is as if $n\{\log(sh(a)) + sg(a)\}$ replace $nsf(a)$. For a fixed $s$, if we can apply the steepest descent's method, the following demonstration will be the same than the precedent where $a \in R^p$.

note $a = (\rho e^{i\vartheta})$, $S(\vartheta, \rho) = \text{Im } \gamma(a) = s \text{ Im } g(a) - I\vartheta + IArgsh(a)$

then, for fixed $s$, we have a uniform code with $\sum_{\square=1}^{\square=p} s_\square \text{ Im } g(a) - \vartheta_\square + IArgsh_\square(a) - \pi k_\square = 0$

$\square = 1, 2, ..., p$. We have $p$ equations making $a$ as a function of the variables $k$, a manifold in $\mathbb{R}^d$ with the uniform code $\kappa_1 : \pi \kappa_1 = s_1 \text{ Im } g_1(a) + IArgsh_1(a) - \vartheta_1$

When $n \to \infty$, we obtain a family of random manifolds depending a uniform random $p$-dimensional vector.

## 5- First analysis critical point $a_c = (a_c, b_c)$

The Perron Frobenius's operator is defined by $\text{Pr } ob(\boldsymbol{a} \in \boldsymbol{B}) = \Sigma_\alpha \text{Pr } ob(f_\alpha^{-1}(\boldsymbol{a} \in \boldsymbol{B}))$, then when we iterate $n$ times, we apply $f^{(n)}$. But the critical points are the images by $f^{(n)}$ of the critical random point $\boldsymbol{a}_c = (a_c, b_c)$.

All the branches are obtained by iterations of the critical point. Thence, the distinct points of cycles $f^{(n)}(a) = a$ will be on these branches. But, heuristically, the iterations of $\boldsymbol{a}_c$ are more and more negligible: iterate $f$ $n$ times, then $\boldsymbol{a}_c$ is iterated $n$ times, but $f^{(k)}(\boldsymbol{a}_c)$ only $PE(n/(k+1))$ times.

## 6- Example: the branches of random parabola of Henon [9]

The iteration of Henon without resonance $f(a, a')$ is characteristic of the phenomena:

Let $\boldsymbol{a} = (a, b)$ and $f(\boldsymbol{a}) = (b + \gamma a - \alpha a^2, \beta a)$. We write $yf(a,b) = x(b + \gamma a - \alpha a^2) + y\beta a$ with $\boldsymbol{y} = (x, y)$. and the resolving deviation for $n = n'$ $e_f^n(\boldsymbol{a}, \boldsymbol{y}) = \partial^n \left[ e^{ay} - e^{yf(a)} \right] / \partial a^n = x^n (y^n - H_n((x\gamma + y\beta)/\sqrt{2\alpha x}))$ where $H_n$ is a Hermitian polynomial. Under the condition $\left| x\gamma + y\beta)/\sqrt{2\alpha x} \right| > |y|$, then, $t = x\gamma + y\beta)/\sqrt{2\alpha x}$ follows the semi circular low. If the initial space we have branches of random parabola following $\beta(1/2, 1/2)$. We must iterate the critical point to obtain the whole random curve. Roughly speaking, in a neighbourhood of the point 0, we obtain random 'parabola' and we can calculate the mean and the variance. In the resonant cases, we must consider the intersection of parabola with the equation of resonance. Rest to examine the decreasing of the branches with the iterations…

## Chapter III – Applications : PDE, ODE, etc.

**I- Introduction, notations and definitions**

Many physical problems can be written with the differential equation $\partial a / \partial t = F(a)$, where $F(a)$ is an application of $R^d$ in $R^d$. In this equation, the unknown variables are a vector $a$ of $R^d$ and the known variables are a vector $t = (q, t_k)$ of $R^{+k}$. The vector $q = (q_1,...,q_i,...q_{k-1}) \in R^{+k-1}$ is the position and $t_k \in R^+$ is the time. The number of equations equals the numbers of unknown variables $a$. We reserve the index $\alpha = 1, 2,...d$ to the unknown $a$ and the index $i = 1, 2,..k$ to the known $t$ with $d \geq k$.

We translate the differential equation as differential iteration to use our methods. As in numerical calculus, we define :

**1- Definition 1**

*A differential iteration is the application $f(a,\delta)$ of $R^d$ in $R^d$, defined by $f(a,\delta) = a + \delta F(a)$ for all fixed $\delta_0 > \delta > 0$ of $R^{+k}$.*

Starting from an initial position $a(0) = (a_\alpha(t(0)))$ at $t(0) = (q(0), 0)$, we iterate $f(a,\delta)$ with a path $\delta$ in all directions, with the particularity that, when $f(a,\delta)$ is iterated $n$ times, $\delta$ will depend on $n$. Let $\delta$ be a small fixed vector of $R^{+k}$. When we iterate $n \in N$ times $f(a,\delta)$, we link $\delta$ to $n: \delta_i = t_i / n$. For fixed $t$, the solution $a(t)$ of the differential equation is obtained when $n \to \infty$. We don't study here the frontier's problems. The main reference on these questions is Arnold [1] and [2].

In order to use the previous results, we define:

**2- Hypothesis H**

*We assume the previous hypothesis of chapters 1 and 2. Quickly said: the polynomial application $F$ iterates a compact set $C \subset R^d$ in $C$ and we can apply the steepest descent's method to $f(a,\delta)$.*

But here, we shall see how the steepest descent's method can be managed.

**3- Difficulties for the probabilistic approach under H**

We want determine the conditions of convergence toward an asymptotic distribution and the distribution induced. But, if we apply all the previous results to $f(a,\delta)$ under H, the distribution will depend on $\delta$. What is the limit of the asymptotic invariant set of zeros $H_n(y,\delta) = \partial^n e^{yf(a,\delta)} / \partial a^n \big|_{a=0}$ when $\delta = t/n$ varies? What happens when $t \to \infty$ ?

We shall proof in the following that, when the Hessian of $F$ is definite, then the invariant set $S(\delta)$ is reduce to fixed points and cycles. The degeneracy of the Hessian can induce a random variable with a complicated behaviour. Here, we study only the case when the iteration is partly linear $f(a,\delta)$.

For ODE, we recall that some deterministic theorems, as the Poincare-Bendixon's, give asymptotic behaviours, but can't be used beyond the dimension 2. One of these theorems says that, if $F$ is Lipschitzian, starting from a point of this compact, we have a local solution, but does not give any information about the behaviour when $t \to \infty$.

**4 - Problem of notation's coherence with definitions**

We shall see a problem when $t \to \infty$ : as the convergence toward some set can depend on the direction $\tau$ of $t$, we must select a direction $\tau = \lim_{t \to \infty} t / |t|$ with $|t| = \Sigma_i t_i$. Hence, $\delta = \tau |t| / n \to 0$ In the case of the ODE or when all the variables $q$ are bounded, then $\tau = 1$.

**5- Partly linear $F$**

**Definition 2**

*The iteration $f(a,\delta)$ is said partly linear if $F$ is partly linear.*

**6- Notation**

As in chapter 2-V-3, we denote: $F(a) = (A(a)b + B(a), C(a)b + D(a))$ where $a = (a,b)$ with $a \in \mathbb{R}^p$ and $b \in \mathbb{R}^{d-p}$. The linear part of $F$ is ($\lambda a$, $\lambda' b$) and all other polynomials have a degree greater than 2. As the phenomenon localizes the distributions around each zero of $F$, we define $\alpha = (\alpha, \beta) \neq 0 = (0,0)$ another zero of $F$ : $F(0) = F(\alpha) = (0,0)$,

- At $0$, we define the PRF of $f(a,\delta)$. Let $G(a) = (a + \tau A(a)b, b)$ be the asymptotic partly linear and $\gamma_G(a)$ be the PRF of $G(a)$ : $\gamma_G(a) = x(a + \tau A(a)b) + yb - n\log a - n\log b$ for $y = (x,y) \in \mathbb{R}^d$.

- Translating $f(a,\delta)$ (proposition 1) the origin at $\alpha$, we note $a = u + \alpha$ with $u = (u,v)$ in $\mathbb{R}^p$ to obtain the partly linear iteration $f_\alpha(u,\delta) = f(u+\alpha,\delta) - \alpha$. Let the asymptotic partly linear: $G_\alpha(u) = (u + \tau A(u + \alpha)v, v)$ at $\alpha$ and $\gamma_{G_\alpha}(u) = x(u + \tau A(u+\alpha)) + y(v - n\log u) - n\log v$ be the PRF of $G_\alpha(u)$.

## II- Theorem 3

*If the differential iteration $f(a,\delta) = a + \delta F(a)$, applying $C \subset \mathbb{R}^d$ in $C$ for all $\delta_0 > \delta > 0$ with $\delta = |t|\tau/n \to 0$ in $\mathbb{R}^{+k}$ is polynomial and apply a compact $C$ in $C$ with two fixed points $0, \alpha \in C$, and verify the hypotheses H of chapter1 and chapter2:*

*1—If $H_n(y,\delta)$ has an dominating asymptotic distribution of zeros, the process converges to this distribution or to fixed point according to $\lambda\tau$ be positive or negative. In the EDO case, the condition becomes $\Sigma\lambda$. In the case of the PDE, the planes $\lambda\tau = 0$ get big gaps.*

*2-- If the Hessian of $yF(a)$ is definite, there exists, besides the fixed points of $F$, only cyclical orbits with uniform density and period $T$ such as $\int_0^T F(a(t+t_0))dt = 0$, $\forall t_0$.*

*3- If $F$ is partly linear, then asymptotically the critical point $a = (a,b)$ de $\gamma_f(a)$ can tend the critical point of either $\gamma_G(a)$ at $0$ or $\gamma_{G_\alpha}$ at $\alpha$. Thence, we can find the asymptotic invariant solution among the solutions of each iteration $G_\alpha(a)$. As $G_\alpha$ are partly linear, we have families of random manifolds by iteration of the critical point.*

*4--According with $\mathrm{Re}(x\tau A(a)\beta)$ is positive or negative at critical point, we shall take the distribution at $0$ or at $\alpha$ that gets get the greatest contribution of the steepest descent's method. The critical points realising $\mathrm{Re}(x\tau A(a)\beta) = 0$ can induce a route of communication between the two sets.*

The proof is guided by the geometry of the attractor of Lorenz. But, as the Hessian is degenerated, we don't obtain a complete solution.

### 1- Proof of point 1: conditions of convergence

The real eigenvalues $\rho$ of the linear part of $f(a,\delta)$ denoted $\lambda \approx 1 + \rho\delta$ verify $|(\partial F(0)/\partial a) - \lambda I| = 0$. As in chapter 1, we write the resolving deviation of the iteration $f(a,\delta)$ : $e_f^n(y,\delta) = y^n - H_n(y,\delta)$ with $H_n(y,\delta) = \partial^n e^{yf(a,\delta)}/\partial a^n |_{a=0}$. For all fixed $\delta > 0$, we can write : $e_f^n(y,\delta) = y^n - \rho^n H_n(y,\delta)/\rho^n \to 0$ when $n \to \infty$. $H_n(y,\delta)/\rho^n$ has its leading term equal 1.

when $n \to \infty$ : $\rho^n = (1 + \lambda\delta)^n \approx e^{n\lambda\delta}$ with $t = \delta = \tau|t|/n$, provided that $|t|/n$ is bounded. That expression tend toward infinite, 0 or 1 according to $\lambda\tau = \Sigma_\square \lambda_\square \tau_\square$ be positive, negative or null. Then, the zeros of $H_n(s,\rho)$ give the limit or not.

In the case of the PDE, consider $\lambda\tau = 0$. The direction $\tau$ can vary and becomes $\lambda\tau = 0$. The process can enter in resonance and conversely. Then, we can observe big gaps in the behaviour.

Note that even though the fixed point is repulsive, it can exist attractive cycles which mask these facts.

### 2 – Proof of point 2: Cycles or fixed points when the Hessian does not degenerate.

- Non-cyclical case; we study now the zeros of $H_{n-1}(y,\delta)$ by the steepest descent's method. We search the critical points of the PRF with : $\gamma(a,\delta) = ya + \delta yf(a,\delta) - n\log a$ where $\delta = \tau|t|/n$. At $0 \in Fix(f)$ the domain $S(\delta)$ is determined by the critical points of $\gamma(a,\delta)$, provided that $\gamma(a,\delta) \to -\infty$ when $(y,n) \to \infty$ and the Hessian be definite negative. The equations of the critical point $a$ are $a_1 \partial \gamma_c,\delta)/\partial a_1 = y_1 a_1 + \delta \partial yF(a)/\partial a_1 - n = 0$, $= 1,2,..,d$. if $\delta \to 0$ and $n \to \infty$, $y_1 a_1 - n \to 0$ and $a_1$ is real. There is no probabilistic invariant distribution.

- Cyclical case:

$f^{(n)}(a,\delta) = a$. We observe that, if $\delta \to 0$ independently of $n$, the cycles disappear.

But, we remark for any cycle $K : f^{(n)}(a,\delta) = f^{(n-1)}(a,\delta) + \delta F \circ f^{(n-1)}(a,\delta) = a$ for all $a \in K$, because we have the limit $f^{(n)}(a,0) = f^{(n-1)}(a,0) = a$ when $n \to \infty$. Nevertheless, for all $\delta = t/n > 0$, we can eventually found solutions with period $T$ such as $\lim_{n \to \infty} f^{(n)}(a,T/n) = a$. In the ODE case, as each point is matched once per period $T$, the asymptotic density is uniform. Moreover, the recurrence's relation can be written: $\partial f^{(n)}(a,\delta)/\partial \delta = [1 + \delta \partial F(u)/\partial u] \partial f^{(n-1)}(a,\delta)/\partial \delta + F \circ f^{(n-1)}(a,\delta)$ en $u = f^{(n-1)}(a,\delta)$. When $\delta \to 0$, we have by recurrence for all $a \in K$ : $(1/n) \sum_{k==0}^{k=n-1} F \square f^{(k)}(a,\delta) \to 0$ and asymptotically $\int_0^T F(a(t+t_0))dt = 0, \forall t_0$. In the PDE case, we obtain similar results.

### 3 - Proof of point3 : Partly linear differential iteration

We come back to the partly linear differential iteration (in chapter 2-V-3): $f(\boldsymbol{a},\delta) = \boldsymbol{a} + \delta F(\boldsymbol{a})$ with $\delta = |t|\tau/n$. The real zeros of $F$ are fixed, isolated points of $f(\boldsymbol{a},\delta)$, and don't depend on $\delta$. We denote $\boldsymbol{a}_1 = (a,b_1)$ with $b_1 = b|t|/n$ and $\boldsymbol{y}_1 = (x,y_1)$ with $y_1|t|/n = y$ such as $y_1 b_1 = yb$

$f(\boldsymbol{a}_1,\delta) = (a + \tau A(a)b_1 + \delta B(a), b_1 n/|t| + \tau C(a)b_1 + \delta D(a))$. Then, with $\boldsymbol{y}_1 = (x,y_1)$:

$\boldsymbol{y}f(\boldsymbol{a},\delta) = x[a + \tau A(a)b_1 + \delta B(a)] + y_1[b_1 + \delta C(a)b_1 + (|t|/n)D(a)]$

As $\delta,|t|/n \to 0$, we obtain asymptotically $\boldsymbol{y}f(\boldsymbol{a},\delta) \to \boldsymbol{y}_1 \boldsymbol{G}(\boldsymbol{a}_1) = xa + \tau A(a)b_1 + y_1 b_1$ and, writing $y,b$ instead of $y_1,b_1$ to lighten the notation, we have the PRF : $\gamma_f(\boldsymbol{a},\delta) \to \boldsymbol{y}\boldsymbol{G}(\boldsymbol{a}) - n\log \boldsymbol{a}$

$n\gamma(\boldsymbol{a}) = n(\boldsymbol{1}\log(sh(\boldsymbol{a})) + sg(\boldsymbol{a}) - \boldsymbol{1}\log \boldsymbol{a})$

This phenomenon localizes the distributions around each zero of $F$ and transforms $\boldsymbol{0}$ in a point fixed of $\boldsymbol{G}(\boldsymbol{a})$. Let $\boldsymbol{\alpha} = (\alpha,\beta)$ another zero of $F$, distinct of $\boldsymbol{0}$ : $F(\boldsymbol{\alpha}) = (0,0)$, the translated iteration (proposition 1 ) of $\boldsymbol{\alpha}$ $f_\alpha(\boldsymbol{a},\delta) = f(\boldsymbol{a} + \boldsymbol{\alpha},\delta) - \boldsymbol{\alpha}$ will be yet partly linear.

### 4 - Proof of point4 : route of communication

Let $\boldsymbol{G}_\alpha(\boldsymbol{a}) = (a + \tau A(a + \alpha)b, b)$ be the asymptotic iteration de $f_\alpha(\boldsymbol{a},\delta)$ localized en $\boldsymbol{\alpha}$. We observe that, for $\delta = 0$, the asymptotic distribution $\boldsymbol{G}_\alpha(\boldsymbol{a})$ en $\boldsymbol{\alpha}$ is not the translated function of $\boldsymbol{G}(\boldsymbol{a})$ en $\boldsymbol{0}$. As we translate $T\boldsymbol{a} = \boldsymbol{a} + \boldsymbol{\alpha}$, $Tf_\alpha T^{-1}$ gives the same invariant measure, then, if $f(\boldsymbol{a},\delta) \to \boldsymbol{G}(\boldsymbol{a})$ and $f_\alpha(\boldsymbol{a},\delta) \to \boldsymbol{G}_\alpha(\boldsymbol{a})$, we should have $T\boldsymbol{G}_\alpha T^{-1} = \boldsymbol{G}$. But, we observe that the iteration computed at $\boldsymbol{\alpha}$ and translated to $\boldsymbol{0}$ : $T\boldsymbol{G}_\alpha T^{-1}(\boldsymbol{a}) = (a + \tau A(a)(b - \beta),b) = \boldsymbol{G}(\boldsymbol{a}) - (\tau A(a)\beta,0)$. That can be explained by the fact that $\boldsymbol{\alpha}$ is not fixed point of $\boldsymbol{G}$.

But, for all $\delta \neq 0$, as $\boldsymbol{\alpha}$ is fixed point of $f(\boldsymbol{a},\delta)$, $f(\boldsymbol{a},\delta)$ and $f_\alpha(\boldsymbol{a},\delta)$ determine a distribution that we have examined in the domination problems: the sign of $\gamma(\boldsymbol{a},\delta) - \gamma_\alpha(\boldsymbol{a},\delta) = x\tau A(a)\beta$ gives the most important contribution.

**Remark:**

The theorem does not give the full asymptotic solution of the differential equations when the Hessian is degenerated. The iteration $G(a)$ can generate cycles and many problems that we have not yet completely solved. We note that the degeneracy of $G(a)$ is very special (one of the worst):

- If we start the iteration from a point $a_0 = (a_0, b_0)$, $G(a)$ will be $G(a) = (a + \tau A(a)b_0, b_0)$. In general, the iteration is without degeneracy if the Hessian of $x\tau A(a)b_0$ is not degenerate. But we don't know $b_0$.

- If we use the steepest descent, the critical point of the PRF $\gamma(a,b) = xa + x\tau A(a)b + yb - n\log ab$ is now defined by two equations: $\partial \gamma(a,b)/\partial a = x + x\tau \partial A(a)\partial ab - n/a = 0$ and $\partial \gamma(a,b)/\partial b = x\tau A(a) + y - n/b = 0$. Then the PRF of $a$ is $\gamma(a,b) = xa - n\log ab$
$n\gamma(a) = n1\log(xA(a) + y) + xa - n1\log a$ where $b$ is now defined by $x\tau A(a) + y - n/b = 0$.

- If we examine our computations, we observe that we suppose implicitly that $b >> a$. Then, the problem splits in two parts: the iteration $a + \tau A(a)b$ and a differential equation $\partial b/\partial t \approx C(a)b$. As we suppose $a \approx 0$, $C(a) \approx C(0)$ and we can write $\partial b/\partial t \approx C(0)b$ or perhaps with a little random noise $\varepsilon$: $\partial b/\partial t \approx C(0)b + \varepsilon$. Then, we can input the asymptotic finite values of the differential equation in the iteration $a + \tau A(a)b$. Does this method lead to the first or the second result?

The steepest descent seems the best…

### III- Examples
#### 1- Deterministic solutions

- **The logistic** $da/dt = \alpha(a - a^2)$. $f(a,\delta) = a + \delta F(a) = a + \delta\alpha(a - a^2)$.

We have 2 fixed points $a = 0$ and $a = 1$. But, only $a = 1$ has $1 - \delta\alpha < 1$. The limit is this point.

- **Hamilton's equations in** $R^d$

- Keeping the notations of the classical mechanics $[1]$:

$dp/dt = -\partial H(p,q)/\partial q$ ; $dq/dt = \partial H(p,q)/\partial p$.

These equations are an application of $R^{2d}$ in $R^{2d}$ when the Hamiltonian $H$ is autonomous. With $F = (-\partial H(p,q)/\partial q, \partial H(p,q)/\partial p)$, we construct the iteration $f(a,\delta) = a + \delta F(a)$, $a = (p,q)$ and assume that it applies a compact C in C for all small $\delta > 0$.

We search the fixed points (of Lagrange or Trojans) $\partial H(p,q)/\partial q = 0$ and $\partial H(p,q)/\partial p = 0$. Then, we compute the eigenvalues of the linear part at these fixed points, denoting $I$ the identity matrix of $R^d$ and vanishing the determinant composed with 4 blocks $d \times d$:

$$J(\lambda) = \begin{vmatrix} -\partial^2 H/\partial q^2 - \lambda I & -\partial^2 H/\partial p\partial q \\ \partial^2 H/\partial p\partial q & \partial^2 H/\partial p^2 - \lambda I \end{vmatrix} = 0$$

Let $\gamma(p,q) = xp - x\delta\partial H(p,q)/\partial q + yq + y\delta\partial H(p,q)/\partial p - n\log p - n\log q$ be the PRF. The invariant density will be seek among imaginary critical points of:

$n/q = y - x\delta\partial^2 H/\partial q^2 + y\delta\partial^2 H/\partial p\partial q$

$n/p = x - x\delta\partial^2 H/\partial q\partial p + y\delta\partial^2 H/\partial p^2$

giving the distribution if the Hessian is degenerated except the case of the cycles with a uniform distribution. Thence $\int_0^T \partial H/\partial q \, dt = 0$ and $\int_0^T \partial H/\partial p \, dt = 0$.

- **k bodies with mass** $m_i$ **in a field of forces derived of a potential** $U(q)$

Hamilton's equations are: $dp/dt = -\partial U/\partial q$, $dq/dt = p/m$ with coordinates $p_i/m_i$ ;

If we can check the validation of the calculus, the PRF of the differential iteration is:
$\gamma(p,q) = xp + yq + \delta(-x\partial U/\partial q + yp/m) - n\log p - n\log q$ which is composed by two functions independent: $\gamma(p,q) = \gamma_0(p) + \gamma_1(q)$ with $\gamma_0(p) = (xp + \delta yp/m - n\log p)$ and $\gamma_1(q) = (yq - \delta x\partial U/\partial q - n\log q)$.

Thence the equations of the critical points $(p,q)$ are:

$\partial \gamma / \partial p = x - \delta x \partial^2 U / \partial p \partial q + \delta y / m - n / p = x + \delta y / m - n / p = 0$

$\partial \gamma / \partial q = y - \delta x \partial^2 U / \partial q^2 - n / q = 0$.

If $\gamma(p,q)$ has a non degenerate Hessian, the system does not have a random solution. After checking calculus ($U(q)$ is not polynomial), if the potential is Newtonian, we have cycles or points.

**2- Probabilistic solutions**

**2-1- Lorenz equations ( see [16] for instance)**

This equation is very important because the partial linearity of its differential iteration induce an asymptotic probabilistic solution. It's the best example to illustrate the previous results.

*- Presentation of the differential iteration and fixed points*

We write these equations with our notations $da/dt = F(a)$, where $a = (a,b,c)$ :
($da/dt = \sigma(b-a)$ ; $db/dt = \rho a - b - ac$ ; $dc/dt = -\beta c + ab$).

With $\delta$ being the path, the differential iteration is quadratic and partly linear in $a$ :
$(a + \delta\sigma(b-a)$ ; $b + \delta(\rho a - b - ac)$ ; $c + \delta(-\beta c + ab))$ The iteration applies a compact $C$ in $C$ for $\delta_0 > \delta > 0$ (the phenomenon happens between a hot plane and a cold plane). The fixed points are the zeros of $F(a) = 0$. When $\rho > 1$ and $\alpha = \sqrt{\beta(\rho-1)}$, we have three: the point $0 = (0,0,0)$, and two symmetrically with the $c$ axis: $\boldsymbol{\alpha}_+ = (\alpha, \alpha, \alpha^2/\beta)$ and $\boldsymbol{\alpha}_- = (-\alpha, -\alpha, \alpha^2/\beta)$.

The eigenvalues equation of the linear part is : $(\beta + \lambda)[(\sigma + \lambda)(1+\lambda) - \sigma\rho] = 0$ at $\boldsymbol{0}$, but $\lambda(\beta+\lambda)(1+\sigma+\lambda) - \alpha^2(\lambda + 2\sigma) = 0$ at $\boldsymbol{\alpha}_+$ or at $\boldsymbol{\alpha}_-$. The coefficients $(\beta,\rho,\sigma)$ are assumed such the three fixed points are repulsive. Hence, we have random distributions around each point fixed. We don't speak here about cycles or resonances. It remains many things to clarify.

**-PRF at the fixed point $\boldsymbol{0}$**

When $\delta \to 0$, the differential iteration tends toward: $G(a) = (a$ ; $b + \rho a - ac$ ; $c + ab)$.
If $y = (x,y,z)$, the PRF is:
$\gamma(a) = yG(a) - nLoga - nLogb - nLogc = yG(a) - nLogabc$

Denoting $y = ns$, with $s = (r,s,t)$ and $\boxed{\varpi = r + s\rho}$,

We can write: $yG(a)/n = a\varpi + bs + ct - sac + tba = L(a) + Q(a)$

with a linear part $L(a) = a\varpi + bs + ct$ and a quadratic part $Q(a) = -sac + tba$

As the application is quadratic, we can compute directly with Hermitian polynomials. The Hessian of $G(a)/n$ has the constant quadric $Q$ and the eigenvectors matrix $T$ is orthogonal for all $a$ :

$Q = \begin{pmatrix} 0 & t & -s \\ t & 0 & 0 \\ -s & 0 & 0 \end{pmatrix}$, with: $T = \dfrac{1}{\mu\sqrt{2}} \begin{pmatrix} 0 & \mu & \mu \\ s\sqrt{2} & -t & t \\ t\sqrt{2} & s & -s \end{pmatrix}$

$\mu$ is the positive eigenvalue of the characteristic equation of $Q$ : $\mu(\mu^2 - s^2 - t^2) = 0$.

**- Changing basis**

The orthogonal application $a = Tu$ avec $u = (u,v,w)$ transforms the 3 factors of the PRF:
$\gamma(a)/n = sG(a) - Logabc = L(a) + Q(a) - Logabc$ :

- $Q(a)$ becomes $-\mu v^2 + \mu w^2$ ,

- $L(a)$ becomes $LTu = \mu u + \varpi v/\sqrt{2} + \varpi w/\sqrt{2}$ .

- As the orthogonal transformation $T$ lets invariant the volume of the rectangular parallelepiped $abc$, we have $Logabc = \log uvw$.

In the new basis $\boldsymbol{u}$, $\gamma(u)/n = \gamma_1(u) + \gamma_2(v) + \gamma_3(w)$ is the sum of three independent functions:
$\gamma_1(u) = \varpi u/\sqrt{2} - \mu u^2 - \log u$ , $\gamma_2(v) = \mu v - \log v$ and $\gamma_3(w) = \varpi w/\sqrt{2} + \mu w^2 - \log w$

becomes in base $u$ : $\gamma(u)/n$ the n derivates of which are with $H_n$ the Hermitian polynomial:

$\partial^n \exp(\varpi u/\sqrt{2} - \mu u^2)/\partial u^n \big|_{u=0} = H_n(\varpi/2\sqrt{\mu})(-\sqrt{2\mu})^n$ ; $\partial^n \exp(\mu v)/\partial v^n \big|_{v=0} = \mu^n$ ;

$\partial^n \exp(\varpi w/\sqrt{2} + \mu w^2)/\partial w^n \big|_{w=0} = H_n(\varpi/2i\sqrt{\mu})(-i\sqrt{2\mu})^n$

Thence, the real zeros of $\partial^n e^{\gamma(u)}/\partial u^n \big|_{u=0}$ are zeros of $H_n(\varpi/2\sqrt{\mu})$.

- **Computing the solution**

As we note : $\partial\gamma(u)/\partial u = \partial\gamma(a)/\partial a . \partial a/\partial u = \partial\gamma(a)/\partial a T$, critical points of $\gamma(a)$ are theses of $\gamma(u)$ and conversely. Hence, it is equivalent to compute the resolving deviation $e_f^n(a,y)\big|_{a=0} = \partial^n(e^{ya} - e^{yf(a,\delta)})/\partial a^n\big|_{a=0}$ in basis $a$ or in basis $u$. The question is now: $\partial^n e^{nsf(a,\delta)}/\partial a^n\big|_{a=0} >> \partial^n(e^{nsa})/\partial a^n\big|_{a=0}$ in the base $u$ ? The point fixed $0$ doesn't be attractive. As $ya/n = sTu = \mu u + (r/\sqrt{2})(v+w)$, the derivates third $n$ of $\exp sTu$ are $\mu^n(r/\sqrt{2})^{2n}$. Thence, the resolving deviation: $\mu^n\left[(r^2/2)^n - H_n(\varpi/2\sqrt{\mu})H_n(\varpi/2i\sqrt{\mu})(i2\mu)^n\right]$

The leading term of the second member of the deviation must be greater than $(r^2/2)^n$ to obtain the probabilistic solution. The condition is $|\varpi| = |r+s\rho| > r$ . Then $\varpi/(2\sqrt[4]{(s^2+t^2)})$ follows asymptotically the Wigner's law. We have a family of random ovals of law beta $\beta(1/2,1/2)$ ,

- **Computation at the other fixed points and dominating distributions**

We search the distributions around the 2 other fixed points. To compute $sG_+(a)$ at $\alpha_+$ ( resp. $sG_-(a)$ at $\alpha_-$ ) from $sG(a)$ at $0$, it is enough to change in the differential iteration $a = (a,b,c)$ in $a + \alpha_+ = (a, b+\alpha, c+\alpha^2/\beta)$ (resp. $a + \alpha_- = (a, b-\alpha, c+\alpha^2/\beta)$ ):
$(a + \delta\sigma(b-a)$ ; $b + \delta(\rho a - b - ac)$ ; $c + \delta(-\beta c + ab)$ ). Then $sG_+(a) = sG(a) - (s(\rho-c)+tb)\alpha$ and $sG_-(a) = sG(a) + (s(\rho-c)+tb)\alpha$ .

The two distributions centred in $\alpha_+$ and $\alpha_-$ dominate the one centred in $0$. $(s(\rho-c)+tb) = 0$ is route of communication between the two others. At $\alpha_+$ we take $\varpi_+ = r + s\rho + s\alpha^2/\beta + t\alpha$ instead of $\varpi$, and $\varpi_- = r + s\rho + s\alpha^2/\beta - t\alpha$ at $\alpha_-$, to obtain the asymptotic geometry.

- **Geometrical interpretation**

The variable random $(r+s\rho)/(2\sqrt[4]{(s^2+t^2)}) = \chi$ gives the distribution of the random ovals at $0, \alpha_+, \alpha_-$. The plane $c = r + s\rho$ intersects the cylinder $(c/2\chi)^4 = s^2 + t^2$ with random radius $(c/2\chi)^2$. Such are the first approach of a mathematical study. To complete, we must compute also the critical points iterates. Thence:

-**Proposition 2**

*An asymptotic solution of the equation of Lorenz is composed by 2 family of random ovals random $(r+s\rho)/(2\sqrt[4]{(s^2+t^2)}) = \chi$ with beta law $\beta(1/2,1/2)$ , centred on the 2 fixed points $\alpha_+$ and $\alpha_-$ , linked by the route of communication $(s(\rho-c)+tb) = 0$ and the distributions of iterates of critical points.*

**2-2 – The Navier Stokes equation of and its random ovals**

This example is the same as the Lorenz's case except the fact that $\delta$ is not uni dimensional. Fefferman [8], writer of the Clay's foundation, feels uncomfortable to formalize mathematically The Navier Stokes equation : prudently, he keeps the concepts of the XIX century where the equations differentials are expended by Fourier's series in order to compact the solution. Moreover, the presentation is trapped by physical considerations, introducing many parameters. Perhaps the physical reality will be as close as possible? For fear of forgetting something, a difficulty will be neglected? We must admit that the mathematicians have matched many difficulties with the dimension 3.

In the Navier Stokes equation, the speed $u$ and the pressure $p$ of a fluid are the unknown variables and the known variables are the position $x \in R^n$ and the time $\tau \in R^+$. The divergence must be null. We have not such considerations here. We take our hypothesis to determine the invariant measure in a compact set, then we suppose $u$ and the pressure $p$ bounded. Note that the spatial periodicity of the solution $a(q,t) = a(q+e,t)$ in $[8]$ induces a space's bound of the recurrence for finite $t$ and $\tau = 1$. Nevertheless, we must translate the equations $[8]$ into iterations, thence obtain $\partial a / \partial t = F(a)$ where the unknown variables are denoted by $a$, and the known variables by $t$. The multiplicity of the variables makes more difficult the translation.

**-Translating $[8]$ into $\partial a / \partial t = F(a)$**

At the beginning, we have as many equations than unknown variables with the divergence null, then, we introduce as many equations than unknown intermediary variables :

- the matrix $\partial u / \partial x = b$ introduces the unknown intermediary $b \in R^n \times R^n$ with an equal number of equations. Thence : $\sum_{j=1}^{j=n} u_j \partial u_i / \partial x_j = ab_i$ with the matrix $U_1$, composed of 0 and 1, we have:

- $div\ u = \sum_{j=1}^{j=n} b_{jj} = U_1 b = 0$

- the field of the external forces $f_i(x,\tau) = f_i(c)$ is here polynomial (or well approximated by polynomials) verifying our hypothesis.

We group together all these unknown intermediary linear variables so various as $\partial u / \partial x = b$, $\Delta u_i = \sum_{j=1}^{j=n} \partial^2 u_i / \partial x_j^2 = \sum_{j=1}^{j=n} \partial b_{ij} / \partial x_j$, $\partial p / \partial x_j = d_j$ or $\partial c / \partial t = 1$ and we denote that in a block: $\partial d / \partial t = B$ where $d$ and $B$ are linear en $(a,b,c)$

- The Navier Stokes equation $[8]$ is written in terms of $\partial a / \partial t = F(a)$ :

$\partial a / \partial \tau = -ab + Ad + f(c)$ where $A$ is a matrix with constant coefficients.

$\partial d / \partial t = B$

**- Solution**

We have obtained a good equation: $\partial a / \partial t = F(a)$ with $a = (a,b,c,d)$ and $t = (x,\tau)$, with a lot of unknown variables linked by linear relations. In the main equations, the only non-linear polynomials are $ab$ and $f(c)$. Thence, the iteration $f(a,\delta) = a + \delta F(a)$ has linearity. Assume that the iteration applies a compact $C$ in $C$ for all $\delta > 0$. With $y = (x,y)$, we can write $P = yF(a,b,c,d)$ : $P = x[-ab + Ad + f(c)] + zc + td + ux$ and $\gamma(a) = P - n\log a$

The Hessian of $P$ is, if ${}^t x$ is the transposed vector of $x$ : $Q = \begin{pmatrix} 0 & {}^t x & 0 \\ x & 0 & 0 \\ 0 & 0 & x\partial^2 f / \partial t^2 \end{pmatrix}$

$x\partial^2 f / \partial t^2$ being computed at 0. The characteristic equation is : $\mu^k(\mu^2 - |x|^2)|x\partial^2 f / \partial c^2 - \mu I| = 0$ where $k$ is such as the total degree of the polynomial equalizes the dimension of the application. The Hessian is degenerated and has a positive eigenvalue $\mu = |x|$) and a negative $\mu = -|x|$ ; the movement is situated in the eigenspaces as in the Lorenz case. In dimension 3, we have an orthogonal matrix as T. We have first groups of solutions with ovals depending on the beta law $\beta(1/2,1/2)$, but we can have independently random sets resulting of the forces $f$. The analysis of this equation is very connected to the Lorenz case: the only new elements are the increasing number of fixed points and the analysis of $f$. However, the more important new things are rifts, which can perturb the behaviour.

**IV- Reflections about the recurrence of Poincaré or the horizon of Lyapounov**
**Definition**

*Let $A$ be a measurable set. A point $a$ is said **recurrent** for A, iif for all integers $p$, exists an integer $k \geq p$ such as: $f^k(a) \in A$.*

1- **Poincare's theorem** says:

*If $f(a)$ is a bijection and preserve the measure μ (for all $A \subset \sup p(\mu) : \mu(A) = \mu(f^{-1}(A)) > 0$ ),*

*Then, for all $A \subset \sup p(\mu)$, the points $a \in A$ are almost recurrent.*

For all small $\delta_0 > \delta > 0$, $f(a,\delta)$ has a reciprocal image. The PF solution is recurrent if exists an invariant measure. As $f(a,\delta)$ applies a compact set $C$ of $R^d$ in C, this measure exists and all the points of the support of the invariant PF distribution are almost recurrent.

Assume that $F(a)$ verify our hypothesis and be partly linear such as the differential iteration be asymptotically equivalent to $G(a) = (a + \tau A(a)b, b)$. The critical points $a$ and their repartition are reciprocal images of the uniform repartition of the unit hypercube $K = (0,1)^p$ compacted with the point 0 with intricacies with the cycles. We divide K in cubes $\kappa_\varepsilon$ with edges $\varepsilon = 1/m$.

Taking a small path $\delta$ such as $f(a,\delta)$ is invertible in the compact. We have an invariant measure.

- We define the variables $t_\delta$ such as the points $a(t_\delta)$ of $f(a,\delta)$ belong at least one time in each $\kappa_\varepsilon$. $t_\delta$ is the variable of first matching for each cube $\kappa_\varepsilon$. For all $t << t_\delta$, the orbit is deterministic.

- If $t_\beta$ is sufficiently large such as we have many points in each cells $\kappa_\varepsilon$, then we can apply the probability calculus and, for all $t >> t_\beta$, the behaviour will be randomized. The ergodism allow us to make calculus as if we draw randomly points in K.

We can define $\theta = (t_\delta, t_\beta)$ as the range between the recurrence of Poincaré and the horizon of Lyapounov.

This conception, nearby the quantum mechanics, gives a transition between determinism and probability. (see Ghys [8]).

**2- Remarks for the mathematics**

Some asymptotic solutions of a differential equation can be randomized. The probabilistic theory brings information where the deterministic approach is can't explain something.

Curiously, the asymptotic distribution is the same, as for ODE as for PDE, but the conditions of convergence are very different. In the PDE case, a little variation of the direction $\tau$ can induce great perturbations in the solutions.

Many ideas on the chaos must be corrected:

- The definition of the chaos, as a situation between determinism and probability, seems devoid of interest. The multiple errors of approximation making unpredictable the behaviour are not a good mathematical argument. The notion of sensibility to the initial conditions is without interest. The notions of local domination or asymptotic behaviours as in Markov process seem better.

- On the other hand, the degeneracy of the Hessian gets a great importance. This degeneracy remains a difficult problem. As example, if we have the equation $da/dt = F(a,t)$ instead of $da/dt = F(a)$, we can add a new coordinate $dt/dt' = 1$ and a new variable $t'$, but that induces linearity with a degeneracy of the Hessian, which modifies the results as we see for the non-autonomous phenomena.

We have a very important lack of mathematical results in this field.

**3- Remarks for the physics**

Limits and bounds are essential for the compactness of the iteration in physics. In many physical problems, the position $q$ in $t = (q,t')$ is bounded. The work of Delabaere and Howls [5] induces these conditions without difficulty for the PRF. In other problems, the spatial periodicity (for example, equation of Navier-Stokes [8]) of the solution $a(q,t') = a(q+e,t')$ induces a space's bound of the recurrence in the space and $\tau = 1$.

Bounds concern the majority of the applications of the physics at finite distance.

**4- Other equations : with decay** $a(\lambda t) = f \circ a(t)$ **or** $f \circ a(t) = \lambda a(t)$

These equations determine frequently the semi-invariant curves of iterations of $f(a)$ where $\lambda$ are the eigenvalues of the linear part of $f$. They often present very curious structures as in the Henon's iteration.

# Chapter 4: Resonance

## I- Introduction

We speak about an important case, we have avoided until now.

It is well known that the resonance makes rigid any mechanical system. It is the same for the solutions of the Perron Frobenius's operator. It is a first approach of this delicate question.

## II- Hypotheses, definition and notations

**Hypothesis: the same as the previous one.**

We diagonalize the linear part of $f$ at the fixed point $0$ : $\partial f(0)/\partial a = \lambda$ real.

**Definition**

*We say that there is resonance if exists an integer vector $k$ such as $\lambda^k = 1$. We call direction of resonance $\zeta_k$, the direction $\zeta_k \log|\lambda| = 0$ in the plane of resonance.*

*We call group of resonance $G^k$ of order $k$, the multiplicative subgroup of such as all $\mu \in G^k$ verify $\mu^k = 1$*

We see that $nk \in \zeta_k$ for $\forall n \in N$, then we take $k$ minimal. As soon as we are in dimension 3, we can have many vectors $k$ which verify such a relation, and we obtain a lattice in the plane of resonance $x \log|\lambda| = 0$. We present here the more simple case where there is one $k$ minimal $\lambda^k = 1$. Thence, we can extend the group $G^k$ on $R^d$ with denoting 1 for the null $k_\square$.

**The problem**

If $f$ apply a compact $C$ in $C$ exists a unique invariant measure, whatever there is resonance or not. We use the resolving deviation $e^n(y) = e^n(y,0) = y^n - H_n(y)$, then, apply the steepest descent to study the zeros of $H_n(y)$. Then we have to compute the set $Z(J_n)$.

When $\lambda^k = 1$, we observe the resolving deviation $e_f^{qk}(y) = y^{qk} - H_{qk}(y)$ where $n = qk$ with $q \in N$. It has the leading term null. We have seen that we can obtain by identification $d$ functions $\phi_q(y) = \sum_{n=0}^{n=kq} b_{n,q} y^n$ vanishing the resolving deviation such as the condition $\phi_q(y) = \phi_{qf}(y)$.

## III -Theorem 4

*Let $k$ be the vector of integers such as $\lambda^k = 1$.*

*The resolving deviation $e_f^{qk}(y)$ remains invariant under the action of the group $G^k$.*

*The point critical $a$ of the PRF verifies the condition of resonance $1 = a^{k-1}$.*

- Let the resolving deviation $e_f^{qk}(y) = y^{qk} - H_{qk}(y)$ when $n = qk$ in the direction of resonance $k$. Let $\mu \in G^k$, and $a = \mu u$, then, as $\mu^{qk} = 1$, we have:

$$\partial^{kq}(e^{y\mu u} - e^{yf(\mu u)})/\partial u^{kq} = \partial^{qk}(e^{ya} - e^{yf(a)})/\partial a^{qk}(\partial a^{kq}/\partial u^{kq}) = \partial^{qk}(e^{ya} - e^{yf(a)})/\partial a^{qk}$$

- The distribution will verify simultaneously $\partial^{qk} e^{yf(a)}/\partial a^{qk} = 0$ for every integer $q$ and $\partial^{n-1} e^{yf(a)}/\partial a^{n-1} = 0$ for all multi indexes such as $n_\square = n$. The PRF for $\partial^{qk} e^{yf(a)}/\partial a^{qk} = 0$ can be written: $\gamma_{qk}(a) = yf(a) - (qk+1)\log a = yf(a) - q(\log a^k + 1/q \log a)$. When $q \to \infty$, $\gamma_{qk}(a) \to sf(a) - \log a^k$ with $s = y/q$. If we compare this result with the PRF for the function $\gamma_n(a) \to \gamma(a) = sf(a) - \log a$, we must take $\log a^k = \log a$ to have the same distribution, and conversely.